\documentclass[AMA,Times1COL]{WileyNJDv5} %STIX1COL,STIX2COL,STIXSMALL

\articletype{Article Type}%

\received{Date Month Year}
\revised{Date Month Year}
\accepted{Date Month Year}
\journal{Journal}
\volume{00}
\copyyear{2023}
\startpage{1}

\usepackage{hyperref}       % hyperlinks
\usepackage{url}            % simple URL typesetting
\usepackage{booktabs}       % professional-quality tables
\usepackage{nicefrac}       % compact symbols for 1/2, etc.
\usepackage{microtype}      % microtypography
\usepackage{lipsum}
\usepackage{graphicx}
\usepackage{cleveref}
\usepackage{amsmath,amsfonts,amssymb}
\usepackage{bbm}
\usepackage[inline]{enumitem}

\newcommand{\oinv}{o_{i}}
\newcommand{\omi}{o_{m}}
\newcommand{\ohyp}{o_{h}}
\newcommand{\ith}{i_{e}}
\newcommand{\hs}{\mathcal{H}}
\newcommand{\diff}{\mathrm{d}}
\newcommand{\iadj}{\mathfrak{a}_i}
\newcommand{\padj}{\mathfrak{a}_p}
\newcommand{\nadj}{\mathfrak{a}_n}

\newcommand{\oadj}{\mathfrak{a}_o}
\newcommand{\obs}{\mathcal{O}}

\newcommand{\Uill}{\mathbf{U}_{R}}

\newcommand{\rec}{\text{rec}}

\usepackage[scientific-notation=true]{siunitx}
\usepackage{mathtools}
\usepackage{subcaption}
\usepackage{caption}
\usepackage{adjustbox}

\begin{document}

\title{Inverse Problem Regularization for 3D Multi-Species Tumor Growth Models}

\author[1]{Ali Ghafouri}

\author[1]{George Biros}

\authormark{Ghafouri \textsc{et al.}}
\titlemark{Inverse Problem Regularization for 3D Multi-Species Tumor Growth Models}

\address[1]{\orgdiv{Oden Institute}, \orgname{University of Texas at Austin}, \orgaddress{\state{TX}, \country{USA}}}

\corres{Ali Ghafouri, 201 E 24th St, Austin, TX 78712, \email{ghafouri@utexas.edu}}

%\presentaddress{This is sample for present address text this is sample for present address text.}

%\fundingInfo{Text}
%\JELinfo{ejlje}

\abstract[Abstract]{We present a multi-species partial differential equation (PDE) model for tumor growth and a an algorithm for calibrating the model from magnetic resonance imaging (MRI) scans. The model is designed for glioblastoma (GBM) brain tumors.
  The modeled species correspond to proliferative, infiltrative, and necrotic tumor cells. The model calibration is formulated as an inverse problem and solved a PDE-constrained optimization method. The data that drives the calibration is derived by a single multi-parametric MRI image. This  a typical clinical scenario for GBMs. 
The unknown parameters that need to be calibrated from data include ten scalar  parameters and the infinite dimensional initial condition (IC) for proliferative tumor cells. This inverse problem is highly ill-posed as we try to calibrate a nonlinear dynamical system from data taken at a single time. To address this ill-posedness, we split the inversion into two stages. First we regularize the IC reconstruction by solving a single-species compressed sensing problem. Then, using the IC reconstruction,  we invert for model parameters using a weighted regularization term. We construct the regularization term by using auxiliary 1D inverse problems.  We apply our proposed scheme to clinical data. We compare our algorithm with single-species reconstruction and unregularized reconstructions. Our scheme enables the stable estimation of non-observable species and quantification of infiltrative tumor cells.}


%\jnlcitation{\cname{%
%\author{Taylor M.},
%\author{Lauritzen P},
%\author{Erath C}, and
%\author{Mittal R}}.
%\ctitle{On simplifying ‘incremental remap’-based transport schemes.} \cjournal{\it J Comput Phys.} \cvol{2021;00(00):1--18}.}

\maketitle

\renewcommand\thefootnote{}
%\footnotetext{\textbf{Abbreviations:} ANA, anti-nuclear antibodies; APC, antigen-presenting cells; IRF, interferon regulatory factor.}

\renewcommand\thefootnote{\fnsymbol{footnote}}
\setcounter{footnote}{1}

\section{Introduction}
\label{sec:intro}
Biophysical models of tumor growth can help identify infiltrative tumorous regions in the brain that are not clearly visible in clinical images. This information can be then used for  patient stratification, preoperative planning, and treatment planning. A main challenge is finding patient-specific model parameters based on limited clinical data, typically just one snapshot before treatment planning. These model parameters include the healthy brain anatomy, the location in which the tumor is started, and growth model coefficients. 
In this paper, we present an inversion methodology that integrates classical inverse problem theory with a multi-species brain tumor growth model. Our approach addresses some challenges of single time snapshot inversion and our numerical results show strengths and weaknesses of the scheme.

\subsection{Contributions}
In our previous study \cite{subramanian2019simulation}, we developed a multi-species tumor growth model integrated with tumor mass deformation but we did not consider its calibration from data.  Calibrating this model is challenging due to the complexity of non-linear coupled partial differential equations (PDEs) and the sparsity of the data. We are interested in reconstruction using a multiparameteric MRI scan at a single time point. 
We split the reconstruction into two stages.  In the first stage we reconstruct the tumor initial condition (IC) assuming a simpler single species model. In the second stage, we reconstruct ten scalar PDE coefficients using the  IC from the fisrt stage. For the first task we use our previous work on inverting for the IC using a single-species PDE \cite{subramanian2020did,scheufele2020fully} and use it to reconstruct model parameters. Our focus in this paper is two fold: (1) analyze and derive algorithms for the second stage inverse prolbem; and (2) evaluate the combined inversion using synthetic data. We summarize our contributions below. 
\begin{enumerate}
  \item Using a 1D model, we analyze the Hessian of the objective function (numerically), a demonstration for the severe ill-posedness of the inversion problem for the 10-parameters multi-species PDE tumor growth model. 
  We also use the 1D formulation to construct a weighted regularization term that we then use for 3D reconstruction. Since we do not have the ground truth parameters for clinical data, we evaluate this scheme with synthetic data. We empirically demonstrate that this regularization term improves model coefficient estimation in the inversion setting. 
  \item We then estimate the tumor's IC using our previous works \cite{subramanian2020did,scheufele2020fully}. We combine IC inversion with model coefficients inversion. We evaluate this full scenario on synthetic data and highlight the strength and weaknesses of the method. 
  \item We parallelize our execution to efficiently solve the forward problem using parallel architecture and GPUs, allowing for 3D inversion on a $160^3$ imaging resolution to take approximately 5 hours for any synthetic case on 4 GPUs.
%Our software will be publicly available upon request. 
%on Github\footnote{\href{https://github.com/ShashankSubramanian/GLIA.git}{https://github.com/ShashankSubramanian/GLIA.git}}. 
\end{enumerate}
As an example, we also report results of the reconstruction using a clinical dataset. A full clinical evaluation of the overall methodology and its comparison with a single-species model is a subject of ongoing investigation and will be reported elsewhere. 

\subsection{Limitations}
The proposed  multi-species includes edema because it is clearly visible in mpMRI scans.  Our edema model is a simple algebraic relation, similar to our previously proposed model in \cite{subramanian2019simulation}. In this analysis, we also ignore the so called mass effect, the normal brain tissue mechanical deformation due to the presence of the tumor. 
We have designed schemes to account for mechanical deformation in single-species models (\cite{subramanian2020multiatlas,subramanian2022ensemble}), but this is beyond the scope of this study. Our growth model does not incorporate diffusion tensor imaging, which can help capture more complex invasion patterns. 
Our single-snapshot multiparameteric MRI does not include more advanced imaging modalities like perfusion. 
In this paper, we use a simplified 1D model to compute the regularization, which differs from the 3D model in some terms. We estimate the tumor IC using a single-species reaction diffusion model proposed in our previous studies \cite{subramanian2020did,scheufele2020fully} and do not invert for the IC using the multi-species model. Our scheme is primarily focused on brain tumors, but it is also applicable to tumors in other organs such as the pancreas, breast, prostate, and kidney.

\subsection{Related works}
The modeling of tumor growth studies can be divided into two categories: 
\begin{enumerate*}
  \item Forward models 
  \item Inverse problems and solutions algorithms
\end{enumerate*}.
Forward models mainly focus on the biophysical model that describes the biological process of brain tumors. Inverse problems and solution algorithms, on the other hand, involve aspects such as noise models, regularization, observation operators, and calibration of forward models in order to find the the parameters to reconstruct data. Despite numerous attempts to simulate tumor growth models, less work has been done on the inverse problems themselves.

The most commonly used forward model for tumor growth is the single-species reaction-diffusion model, as seen in various studies \cite{hogea2007modeling,hogea2008brain,jbabdi2005simulation,mang2012biophysical}. Other models by incorporating mass effect, chemotaxis, and angiogenesis are also explored in \cite{hormuth2018mechanically,saut2014multilayer,swanson2011quantifying}. Despite offering insight into biological processes, these models lack calibration due to mathematical and computational complexities. Our group has worked to address these challenges, enabling clinical analysis of the models \cite{subramanian2022ensemble,gholami2019novel}.

Inverse problems have been extensively studied to calibrate single-species reaction-diffusion models and quantify infiltrative tumor cells \cite{scheufele2020image,fathi2018characterization,peeken2019deep}. In reference \cite{gholami2016inverse}, the inverse problem is formulated by including anisotropy diffusion derived from diffusion tensor imaging. In reference \cite{lipkova2019personalized} the authors utilize bayesian inversion for personalized inversion of model parameters for invasive brain tumors. Our previous works \cite{subramanian2020did,scheufele2020fully} addressed the challenges of calibrating the model coefficients and IC inversion simultaneously, proposing a $\ell_0$ constraint to localize IC inversion. In reference \cite{subramanian2020multiatlas}, the model is combined with mass deformation and clinical assessment is performed for a large number of MRI scans~\cite{subramanian2022ensemble}, improving survival prediction for patients. However, the single-species model does not use all available information from MRI scans and considers necrotic and enhancing cells as single-species while ignoring the edema cells, ignoring the biological process between species and resulting in a model output that differs from observed segmented tumor from MRI scans. 

Recently~\cite{hormuth2021image}, a dual-species model using longitudinal MRI scans to calibrate the model coefficients has been developed. However, longitudinal data is rarely available in clinical practice. In our previous work \cite{subramanian2019simulation}, we formulate a multi-species model from \cite{saut2014multilayer} and combine it with mass deformation of the brain as a forward solver. The solver developed in this study was applied to analyze multiple clinical cases, as detailed in \cite{ghafouri20233d}. The present paper, however, focuses on addressing the inverse problem, with particular emphasis on examining its ill-posed nature. Furthermore, we propose a regularization scheme to mitigate the challenges associated with ill-posedness, an aspect not explored in our previous work \cite{ghafouri20233d}. This investigation aims to enhance the robustness and reliability of the solver in clinical applications.

\subsection{Outline}
In \Cref{sec:methods}, we outline the forward problem and mathematical formulation of the inverse problem. We then present the 1D model, derive the adjoint and gradient formulations and estimate the Hessian of the problem. We then propose an algorithm to compute a weighted regularization operator, which we the use in 3D reconstructions. The solution algorithm for the 3D inversion problem is also discussed, including the estimation of IC and model coefficients in 3D. 
In \Cref{sec:results}, we analyze the ill-posedness of the 1D problem and present inversion results using our computed regularization. The inversion in 3D is also analyzed for two cases: (1) estimating tumor growth parameters given IC and brain anatomy, and (2) estimating both the brain anatomy and IC while inverting for the tumor growth parameters. Results of the full inversion algorithm are presented on two clinical cases, with a comparison of the reconstructed observed species.

\section{Methodology}
\label{sec:methods}
\label{sec:methods}
In the single-species model, we use the aggregate variable $c(\mathbf{x},t)$ to represent all the abnormal tumor cells. Here $\mathbf{x}$ is a point 
in $\mathbb{R}^{d}$ ($d=3$ for 3D model and $d=1$ for 1D model) and $t \in (0, T]$ is the time. We denote $T$ as the time-horizon (the time since the tumor growth onset) and we set $T=1$ to non-dimensionalize PDE forward model for the single-snapshot inversion. In the multi-species model, we have three species: necrotic $n(\mathbf{x},t)$, proliferative $p(\mathbf{x},t)$, infiltrative $i(\mathbf{x},t)$. We can relate the single-species model to the multi-species model by setting $c(\mathbf{x},t) = n(\mathbf{x},t) + p(\mathbf{x},t) + i(\mathbf{x},t)$. 
We consider probability maps related to MRI scans instead of precise tissue types in our model. The brain is described by separate maps for gray matter ($g(\mathbf{x},t)$), white matter ($w(\mathbf{x}, t)$), and cerebrospinal fluid (CSF)-ventricles (VT), denoted by $f(\mathbf{x},t)$. We use diffeomorphic image registration to segment the healthy brain into gray matter, white matter, CSF and VT \cite{gholami2019novel}.
\subsection{Formulation (Forward problem)}
\begin{figure}
  \centering  
  \includegraphics[width=0.9\textwidth]{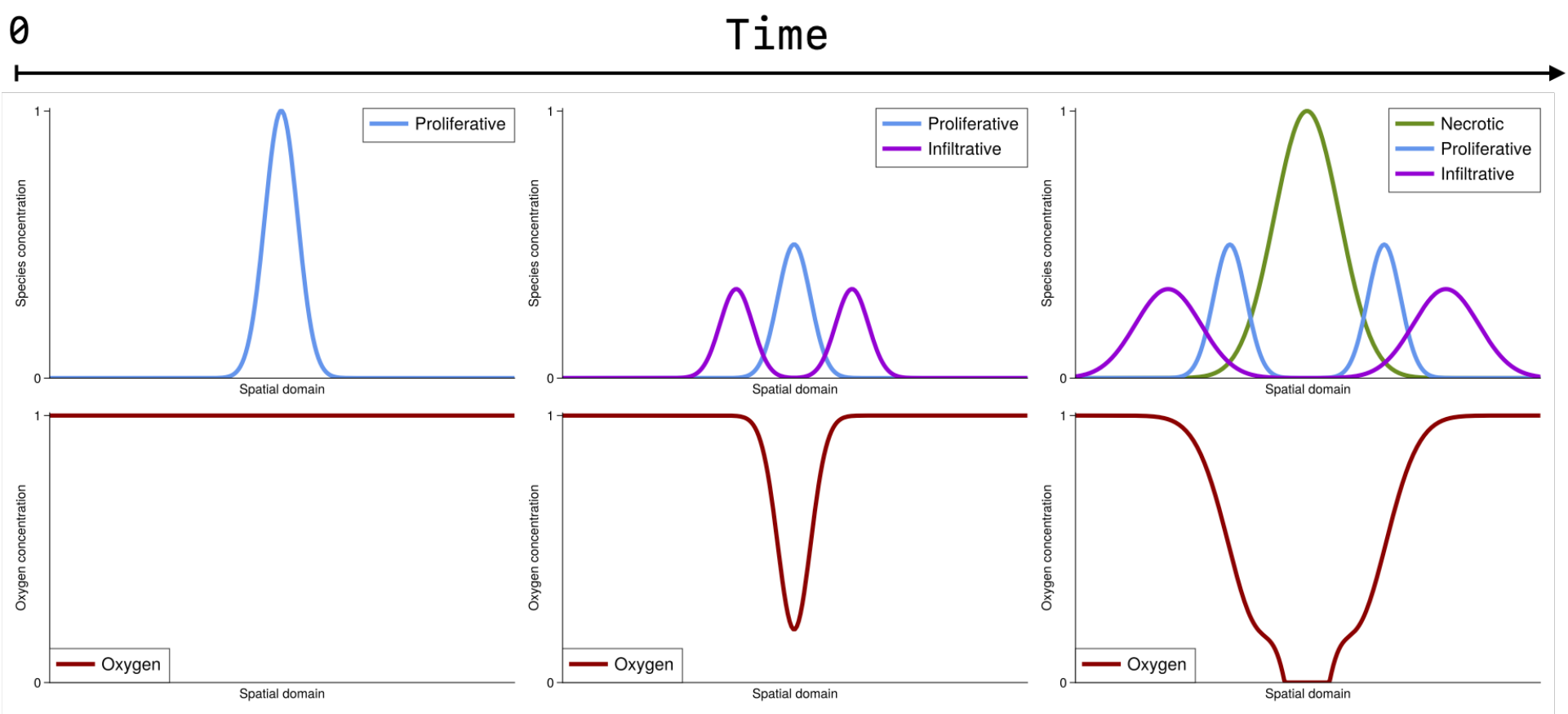}
  \caption{A schematic representation of the multi-species model dynamics in 1D. The tumor species growth process is represented, with vascularization indicated by oxygen.}
  \label{fig:model}
\end{figure}

Our model is a go-or-grow multi-species tumor model that accounts for tumor vascularization. Tumor cells are assumed to exist in one of two states: proliferative or invasive. In an environment with sufficient oxygen, tumor cells undergo rapid mitosis by consuming oxygen. However, under low oxygen concentrations (hypoxia), the cells switch from proliferative to invasive. The invasive cells migrate to richer higher-oxygen areas and then switch back to the proliferative state. 
When the oxygen levels drop below a certain level, the tumor cells become necrotic. We show the schematic progression of multi-species process in \Cref{fig:model} (presented in 1D). Initially, only proliferative cells are present, with the oxygen concentration at its maximum initial value of one. Subsequently, the proliferative cells consume the oxygen, causing them to transition into infiltrative cells. These infiltrative cells then migrate in search of additional oxygen. 
Our model is summarized as,
\begin{subequations}
  \label{eq:forward_problem}
\begin{align}
  \partial_t p - \mathcal{R} p + \alpha(o) p (1 - i) - \beta(o) i (1 - p) + \gamma(o) p (1 - n) & = 0  \quad \text{ in } \Omega \times (0, 1] \label{eq:p_pde} \\ 
  \partial_t i - \mathcal{D} i - \mathcal{R} i - \alpha(o) p ( 1- i) + \beta(o) i (1 - p) + \gamma(o) i ( 1- n) & = 0 \quad \text{ in }  \Omega \times (0, 1] \label{eq:i_pde} \\
  \partial_t n - \gamma(o) (i + p) (1 - n) & = 0  \quad \text{ in } \Omega \times (0, 1] \label{eq:n_pde}\\ 
  \partial_t o + \delta_c o p - \delta_s (1 - o) (w + g) & = 0 \quad \text{ in } \Omega \times (0, 1] \label{eq:o_pde} \\ 
  \partial_t w + \frac{w}{w + g} ( \mathcal{D}i + \mathcal{R}p + \mathcal{R}i) & = 0  \quad \text{ in } \Omega \times (0, 1] \label{eq:w_pde} \\ 
  \partial_t g + \frac{g}{w + g} ( \mathcal{D}i + \mathcal{R}p + \mathcal{R}i) & = 0 \quad \text{ in } \Omega \times (0, 1] \label{eq:g_pde} \\ 
  p(\mathbf{x}, 0) - p_0 & = 0  \quad \text{ in } \Omega \label{eq:p_init} \\ 
  i(\mathbf{x}, 0) & = 0  \quad \text{ in } \Omega \\ 
  n(\mathbf{x}, 0) & = 0  \quad \text{ in } \Omega \\
  o(\mathbf{x}, 0) - 1 & = 0 \quad \text{ in } \Omega \\ 
  g(\mathbf{x}, 0) - g_0 & = 0 \quad \text{ in } \Omega \\
  w(\mathbf{x}, 0) - w_0 & = 0  \quad \text{ in } \Omega \\
  f(\mathbf{x}, 0) - f_0 & = 0   \quad \text{ in } \Omega 
\end{align}
\end{subequations}
The paper uses notations listed in \Cref{tab:notation}. The model ignores mass effects from proliferative and necrotic cells and uses conservation equations with transitional ($\gamma, \beta$ and $\alpha$), source ($\delta_c$, $\delta_s$ and $\mathcal{R}$), and diffusion ($\mathcal{D}$) terms. See \cite{subramanian2019simulation,saut2014multilayer} for model choice and description. The reaction operator ($\mathcal{R}$) for species $s$ is modeled as an inhomogeneous nonlinear logistic growth operator as,
\begin{align}
  \mathcal{R}s = \begin{cases}
    \rho_m s (1 - c),  & o > \oinv \\
    (\cfrac{o - \omi}{\oinv - \omi})\rho_m s ( 1 - c), & \oinv \geq o \geq \omi \\ 
    0, & o < \omi
  \end{cases}
  \label{eq:reac}
\end{align}
where $\omi$ and $\oinv$ are the mitosis and invasive oxygen thresholds, respectively. We use $\omi = \frac{\oinv + \ohyp}{2}$ from \cite{subramanian2019simulation,saut2014multilayer}.
$\rho_m$ is an inhomogeneous spatial growth term defined as,
\begin{align}
  \rho_m (\mathbf{x}) = \rho w + \rho_g g,
\end{align}
where $\rho$ and $\rho_g$ are scalar coefficient quantifying proliferation in white matter and gray matter, respectively. 
We assume $20\%$ growth in gray matter ($\rho_g = 0.2 \rho $) compared to white matter \cite{lipkova2019personalized, swanson2000quantitative, gooya2012glistr}. The growth terms converge to zero when the total tumor concentration ($c$) reaches $1$. The inhomogeneous isotropic diffusion operator $\mathcal{D}$ is defined as follows.
\begin{align}
  \mathcal{D} i = \text{div}(k \nabla i) \label{eq:diff}
\end{align} 
where $k$ defines the inhomogeneous diffusion rate in gray and white matter and we define it as,
\begin{align}
  k = \kappa w + \kappa_g g
\end{align} 
where $\kappa$ and $\kappa_g$ are scalar coefficient quantifying diffusion in white matter and gray matter, respectively. 
We set $\kappa_g = 0.2 \kappa$ \cite{lipkova2019personalized, swanson2000quantitative, gooya2012glistr}.
Note that the migration and proliferation occur only into white/gray matter and not in CSF/VT. The parameters $\alpha, \beta$, and $\gamma$ functions between are set as follows,
\begin{align}
  \alpha(o) & = \alpha_0 \hs(\oinv - o), \label{eq:alpha} \\ 
  \beta(o) & = \beta_0 \hs(o - \oinv), \label{eq:beta} \\ 
  \gamma(o) & = \gamma_0 \hs(\ohyp - o), \label{eq:gamma}
\end{align} 
where the scalar coefficients controlling species transitioning are $\alpha_0, \beta_0$, and $\gamma_0$, and $\hs$ is the smoothed Heaviside function.
\begin{table}[!ht]
  \centering
  \caption{Common notations used to describe the forward and inverse problem for the multi-species tumor growth model.}
  %\resizebox{0.5\textwidth}{!}{
      \begin{tabular}{cll}
          \toprule
          Notation & Description & Range \\           
          \midrule
          $\mathbf{x}$ & Spatial coordinates & ($[0, 2\pi]^3$ for 3D and $[0, 2\pi]$ in 1D)\\ 
          $t$ & Time & $[0, 1]$ \\           
          $p$ & Proliferative tumor cell & $p(\mathbf{x},t) \in [0, 1]$ \\ 
          $i$ & Infiltrative tumor cell  & $i(\mathbf{x},t) \in [0, 1]$ \\ 
          $n$ & Necrotic tumor cell  & $n(\mathbf{x},t) \in [0, 1]$ \\ 
          $c$ & Tumor cell  ($c = i + p + n$) & $c(\mathbf{x},t) \in [0, 1]$ \\ 
          $o$ & Oxygen concentration & $o(\mathbf{x},t) \in [0, 1]$ \\ 
          $g$ & Gray matter cells & $g(\mathbf{x},t) \in [0, 1]$ \\ 
          $w$ & White matter cells & $w(\mathbf{x},t) \in [0, 1]$ \\ 
          $f$ & CSF/VT cells & $f(\mathbf{x},t) \in [0, 1]$ \\ 
          $p_0$ & Initial condition (IC) for proliferative cells ($p(\mathbf{x}, 0)$) (See \Cref{eq:p_init}) & $p(\mathbf{x},0) \in [0, 1]$ \\ 
          \midrule
          $\mathcal{R}$ & Tumor growth operator (logistic growth operator) & - \\ 
          $\mathcal{D}$ & Tumor migration operator (diffusion)  & - \\ 
          $\rho_m$ & Inhomogeneous reaction rate & - \\
          $k$ & Inhomogeneous diffusion rate & - \\
          $\alpha(o)$ & Transition function from $p$ to $i$ cells & - \\ 
          $\beta(o)$ & Transition function from $i$ to $p$ cells & - \\ 
          $\gamma(o)$ & Transition function from $i$ and $p$ to $n$ cells & - \\ 
          \midrule 
          $\rho$ &  Reaction coefficient (See \Cref{eq:reac}) & $[\num{5.0}, \num{25.0}]$ \\ 
          $\kappa$ &  Diffusion coefficient (See \Cref{eq:diff}) & $[\num{0.001}, \num{0.1}]$\\ 
          $\alpha_0$ &  Transitioning coefficient from $p$ to $i$ cells (See \Cref{eq:alpha}) & $[\num{0.1}, \num{10.0}]$ \\
          $\beta_0$ &  Transitioning coefficient from $i$ to $p$ cells (See \Cref{eq:beta}) & $[\num{0.1}, \num{15.0}]$ \\
          $\gamma_0$ &  Transitioning coefficient from $i$ and $p$ to $n$ cells (See \Cref{eq:gamma}) & $[\num{1.0}, \num{20.0}]$ \\
          $\delta_c$ & Oxygen consumption rate (See \Cref{eq:o_pde}) & $[\num{1.0}, \num{20.0}]$  \\ 
          $\delta_s$ & Oxygen supply rate (See \Cref{eq:o_pde}) & $[\num{1.0}, \num{8.0}]$ \\ 
          $\ohyp$ & Hypoxic oxygen threshold (See \Cref{eq:reac,eq:gamma}) & $[\num{0.001}, \num{0.8}]$ \\ 
          $\oinv$ & Invasive oxygen threshold (See \Cref{eq:alpha,eq:beta,eq:reac}) & $[\num{0.2}, \num{1.0}]$ \\
          $\ith$ &  Threshold coefficient for edema (See \Cref{eq:obs_l}) & $[\num{0.001}, \num{0.3}]$ \\ 
          \hline
      \end{tabular}
  %}
  \label{tab:notation}
\end{table}
Note that our model differs from the one in \cite{subramanian2019simulation} in the following ways:
\begin{enumerate}[label=(\roman*)]
\item The transition term $\alpha p $ is changed to ${\alpha} p (1 - i)$ to prevent transitioning once $i \approx 1$ (See \Cref{eq:i_pde,eq:p_pde}).
\item The transition term ${\beta} i$ is changed to ${\beta} i (1 - p)$ to prevent transitioning once $p \approx 1$ (See \Cref{eq:i_pde,eq:p_pde}).
\item The necrosis term ${\gamma} i$ is changed to ${\gamma} i ( 1- n)$ to prevent death once $ n \approx 1$ (See \Cref{eq:i_pde,eq:p_pde,eq:n_pde}).
\item A term in reaction operator $\mathcal{R}$ on a species $s$ is changed from $(1-s)$ to $(1-c)$, where $c$ is the total tumor concentration $(c=i+p+n)$. This change ensures that there will be no reaction or increase of tumor concentration once the total tumor concentration $c \approx 1$ (See \Cref{eq:i_pde,eq:p_pde,eq:w_pde,eq:g_pde}). 
\item In \Cref{eq:n_pde}, the term $\gamma (g+w)$ has been removed to ensure that only infiltrative and proliferative cells can transition into the necrotic state. 
\item We modify the $\beta$ function to be the complement of the $\alpha$ function (See \Cref{eq:beta}) in oxygen terms.
\item We have incorporated the edema model into the observation operator instead of using a separate model (See \Cref{eq:obs_l}).
\end{enumerate}
\subsection{Observation Operators}
We now discuss how we relate problem predictions to MRI images. 
Currently, MRI cannot fully reveal tumor infiltration. Instead, we use a relatively standard preprocessing workflow \cite{bakas2018identifying} in which the MRI is segmented to different tissue types. 
To reconcile the species concentrations generated by our multi-species model with the observed species segmentation obtained from MRI scans, 
we need to introduce additional modeling assumptions related to the so-called \emph{observation operator}. 
These operators map the model-predicted fields to a binary (per tissue type) segmentation, which can then be compared to the MRI segmentation. 
In particular, we define observation operators for the total tumor region ($c$) and each individual species, allowing us to determine the most probable species at each point $\mathbf{x}$. The observation operators are defined as follows:
\begin{align}
  \obs c = \hs(c - w) \hs(c - g) \hs(c - f) \label{eq:obs_c} \\ 
  \obs p = \hs(p - n) \hs(p - i) \obs c \label{eq:obs_p} \\ 
  \obs n = \hs(n - p) \hs(n - i) \obs c \label{eq:obs_n} \\ 
  \obs l = (1 - \obs p - \obs n) \hs(i - \ith) \label{eq:obs_l}
\end{align}
where $\obs c, \obs p, \obs n$ and $\obs l$ are the observed tumor ($c$), proliferative ($p$), necrotic ($n$) and edema cells ($l$), respectively. Note that the species used the above setup are the species at $t=1$. In this setup, we first determine the most probable location for the tumor ($c$) using \Cref{eq:obs_c}, and then we determine the most probable species within $\obs c$ using \Cref{eq:obs_p,eq:obs_n}. To model edema, we adopt a simplified approach similar to \cite{subramanian2019simulation}. In this model, locations with infiltrative concentration above a threshold ($\ith$) are considered as edema if they are not detected as necrotic or proliferative. This model treats edema as a label that is independent of the other PDEs, and can thus be computed in the observation operator $\obs l$ using \Cref{eq:obs_l}. 
To model the Heaviside function $\hs$, we use a sigmoid approximation function given by,
\begin{align}
\hs(x) = \frac{1}{1 + e^{-\omega x}}
\end{align}
where $\omega$ is the shape factor, determining the smoothness of the approximation. This parameter is set based on the spatial discretization. Thus, the proposed observation operators do not have any free parameters. Next, we discuss the inverse problem to invert for the model parameters.
\subsection{Inverse Problem}
Our model requires three parameters to characterize tumor progression:
\begin{enumerate}[label=(\roman*)]
  \item Healthy precancerous brain segmentation ($\Omega$) 
  \item IC of proliferative tumor cells ($p_0$)
  \item Model coefficients vector denoted by $\mathbf{q}$ consisting of the following terms : diffusion coefficient $\kappa$ in \Cref{eq:diff}, reaction coefficient $\rho$ in \Cref{eq:reac}, transition coefficients $\alpha_0$, $\beta_0$ and $\gamma_0$ in \Cref{eq:alpha,eq:beta,eq:gamma}, oxygen consumption rate $\delta_c$ in \Cref{eq:o_pde}, oxygen supply rate $\delta_s$ in \Cref{eq:o_pde}, hypoxic oxygen threshold $\ohyp$ in \Cref{eq:gamma,eq:reac}, 
  invasive oxygen threshold in \Cref{eq:reac,eq:beta,eq:alpha} and threshold coefficient for edema $\ith$ in \Cref{eq:obs_l}. 
\end{enumerate}
We discuss in \Cref{sssec:brain_anatomy} and \Cref{sssec:IC} how we estimate $\Omega$ and $p_0$. Here, we focus on model coefficients $\mathbf{q}$ inversion. Given the $\Omega$ and $p_0$, we frame the inverse problem as a constrained optimization problem,
\begin{align}
  \min_{\mathbf{q}} \mathcal{J} \coloneqq \frac{1}{2} \| \obs p(\mathbf{x},1) - p_{\text{d}} \|^2_{L_2(\Omega)} & + \frac{1}{2} \| \obs n(\mathbf{x},1) - n_{\text{d}} \|^2_{L_2(\Omega)} 
  \frac{1}{2} \| \obs l(\mathbf{x},1) - l_{\text{d}} \|^2_{L_2(\Omega)} \label{eq:inv_problem}
  \\ 
  \text{ s.t. } & \mathcal{F}(p_0,\mathbf{q}) \quad \text{ in } \Omega \times (0, 1] \label{eq:inv_cons}
\end{align}
We aim to minimize the objective function $\mathcal{J}$ by optimizing the scalar coefficients vector $\mathbf{q}$ in the forward model $\mathcal{F}(p_0, \mathbf{q})$. This inverse problem is a generalization of the single-species reaction-diffusion PDE model that has already been shown to be exponentially ill-posed \cite{subramanian2020did,ozisik2018inverse,cheng2008quasi,zheng2014recovering}. We circumvent the ill-posedness by designing a problem regularization that involved a combination of sparsity constraints, a $\max$ constraint in the initial condition and $\ell_2$ penalty on the unknown parameters, $\kappa$ and $\rho$.
As we have more undetermined parameters, we cannot use the regularization scheme in \cite{subramanian2020did,scheufele2020fully}. 
To address this ill-posedness, we first perform analysis on Hessian of the optimization problem and then generalize it to the full 3D inversion model. Next, we perform an analysis for the 1D inverse problem and compute a regularization term to assist in estimation of model coefficients $\mathbf{q}$. 
\subsection{Coefficients inversion analysis for a 1D model}
To demonstrate the ill-posedness of our model, we analyze a multi-species model in 1D. To simplify the problem, we assume that the only normal tissue type is white matter (removing the gray matter, ventricles (VT), and CSF from the model). In this setting, the forward model with periodic boundary conditions is described by:
\begin{subequations}
  \label{eq:forward_problem_1D}
\begin{align}
  \partial_t p - \mathcal{R} p + \alpha p(1 - i) - \beta i ( 1 - p) + \gamma p (1 - n) & = 0 \quad \text{ in } \Omega \times (0, 1]\\ 
  \partial_t i - \mathcal{D} i - \mathcal{R} i - \alpha p(1-i) + \beta i ( 1 - p) + \gamma i (1 - n) & = 0 \quad \text{ in } \Omega \times (0, 1]\\ 
  \partial_t n - \gamma (i + p) (1 - n) & = 0 \quad \text{ in } \Omega \times (0, 1]\\ 
  \partial_t o + \delta_c o p - \delta_s (1 - o) w &= 0 \quad \text{ in } \Omega \times (0, 1]\\ 
  \partial_t w + \mathcal{D} i + \mathcal{R} p + \mathcal{R} i &= 0 \quad \text{ in } \Omega \times (0, 1] \\
  p(x,0) - p_0 &= 0 \quad \text{ in } \Omega \\ 
  i(x,0) &= 0  \quad \text{ in } \Omega \\
  n(x,0) &= 0 \quad \text{ in } \Omega \\
  o(x,0) &= 1 \quad \text{ in } \Omega \\
  w(x,0) - 1 + p_0 &= 0 \quad \text{ in } \Omega
\end{align}
\end{subequations}
The optimization problem has exactly the same form as \Cref{eq:inv_problem,eq:inv_cons},
\begin{align}
  \min_{\mathbf{q}} \mathcal{J} \coloneqq \frac{1}{2} \| \obs p(\mathbf{x},1) - p_{\text{d}} \|^2_{L_2(\Omega)}  + \frac{1}{2} \| \obs &n(\mathbf{x},1) - n_{\text{d}} \|^2_{L_2(\Omega)} 
  \frac{1}{2} \| \obs l(\mathbf{x},1) - l_{\text{d}} \|^2_{L_2(\Omega)} \\ 
  \text{ s.t. }  \mathcal{F}(p_0,\mathbf{q}) & \quad \text{ in } \Omega \times (0, 1]
\end{align}
where $p_{\text{d}}, n_{\text{d}}$ and $l_{\text{d}}$ are the proliferative, necrotic and edema data. The corresponding Lagrangian of this problem is given by :
\begin{align}
  \mathcal{L}  = \mathcal{J} 
  & + \int_{0}^{1} \int_{\Omega} \mathfrak{a}_p [\partial_t p - \mathcal{R} p + \alpha p(1 - i) - \beta i ( 1 - p) + \gamma p (1 - n)] ~\diff \Omega ~\diff t  \\ 
  & + \int_{0}^{1} \int_{\Omega} \mathfrak{a}_i [\partial_t i - \mathcal{D} i - \mathcal{R} i - \alpha p(1-i) + \beta i ( 1 - p) + \gamma i (1 - n)] ~\diff \Omega 
  ~\diff t \nonumber \\ 
  & + \int_{0}^{1} \int_{\Omega} \mathfrak{a}_n [\partial_t n - \gamma (i + p) (1 - n) ] ~\diff \Omega ~\diff t \nonumber \\ 
  & +  \int_{0}^{1} \int_{\Omega} \mathfrak{a}_o [\partial_t o + \delta_c o p - \delta_s (1 - o) w ] ~\diff \Omega ~\diff t \nonumber \\ 
  & + \int_{0}^{1} \int_{\Omega} \mathfrak{a}_w [\partial_t w + \mathcal{D} i + \mathcal{R} p + \mathcal{R} i ] ~\diff \Omega ~\diff t \nonumber \\ 
  & + \int_{\Omega} \mathfrak{a}_p(x,0) (p(x,0) - p_0) ~\diff \Omega \nonumber \\ 
  & + \int_{\Omega} \mathfrak{a}_i(x,0) (i(x,0)) ~\diff \Omega \nonumber \\ 
  & + \int_{\Omega} \mathfrak{a}_n(x,0) (n(x,0)) ~\diff \Omega \nonumber \\ 
  & + \int_{\Omega} \mathfrak{a}_o(x,0) (o(x,0)-1) ~\diff \Omega \nonumber \\ 
  & + \int_{\Omega} \mathfrak{a}_w(x,0) (w(0) - 1 + p_0) ~\diff \Omega \nonumber
\end{align}
where $\mathfrak{a}_p, \mathfrak{a}_i, \mathfrak{a}_n, \mathfrak{a}_o$ and  $\mathfrak{a}_w$ are the adjoint variables with respect to $p, i, n, o$ and $w$. \newpage
The first order optimality conditions can be derived by requiring stationary of the Lagrangian with respect to the state variables ($p,i,n,o$ and $w$). Therefore, we obtain, 
\begin{subequations}
\begin{align}
  \frac{\partial \mathcal{L}}{\partial p} = 0 \\
  \frac{\partial \mathcal{L}}{\partial i} = 0 \\
  \frac{\partial \mathcal{L}}{\partial n} = 0 \\
  \frac{\partial \mathcal{L}}{\partial o} = 0 \\
  \frac{\partial \mathcal{L}}{\partial w} = 0
\end{align}
\end{subequations}
Therefore, the adjoint equations are 
\begin{subequations}
  \label{eq:adjoint_equations_1D}
\begin{align}
  - \partial_t \mathfrak{a}_p  
  - \frac{\partial \mathcal{R}p}{\partial p} (\mathfrak{a}_p - \mathfrak{a}_w) 
  + \alpha (1- i) (\mathfrak{a}_p - \mathfrak{a}_i) 
  + \beta i (\mathfrak{a}_p - \mathfrak{a}_i) 
  + &\gamma (1 - n) (\mathfrak{a}_p - \mathfrak{a}_n) \nonumber \\  
  - \frac{\partial \mathcal{R}i}{\partial p} (\mathfrak{a}_i - \mathfrak{a}_w)   + \delta_c o \mathfrak{a}_o & = 0 
    \text{ in }\Omega \times (0,1],
  \\ 
   \mathfrak{a}_p(\mathbf{x}, 1) + (\obs p - p_{\text{d}}) \frac{\partial \obs p}{\partial p} 
   + (\obs n - n_{\text{d}}) \frac{\partial \obs n}{\partial p} 
   + (\obs l - l_{\text{d}}) \frac{\partial \obs l}{\partial p} & = 0 
     \text{ in }\Omega,
  \\ 
  - \partial_t \mathfrak{a}_i - \mathcal{D} \mathfrak{a}_i + \mathcal{D} \mathfrak{a}_w 
  + \alpha p (\mathfrak{a}_i - \mathfrak{a}_p) 
   + \beta (1 - p) (\mathfrak{a}_i - \mathfrak{a}_p) 
  + \gamma (1 &- n) (\mathfrak{a}_i - \mathfrak{a}_n) \nonumber\\ 
  - \frac{\partial \mathcal{R} i}{\partial i} (\mathfrak{a}_i - \mathfrak{a}_w) 
  - \frac{\partial \mathcal{R} p}{\partial i} (\mathfrak{a}_p - \mathfrak{a}_w) & = 0 
    \text{ in }\Omega \times (0,1],
   \\ 
   \mathfrak{a}_i(\mathbf{x}, 1) + (\obs p - p_{\text{d}}) \frac{\partial \obs p}{\partial i} 
   + (\obs n - n_{\text{d}}) \frac{\partial \obs n}{\partial i} 
   + (\obs l - l_{\text{d}})  \frac{\partial \obs l}{\partial i} &= 0 
     \text{ in }\Omega 
  \\ 
   -\partial_t \mathfrak{a}_n 
  + \gamma ((i+p) \mathfrak{a}_n - i \mathfrak{a}_i - p \mathfrak{a}_p) 
  - \frac{\partial \mathcal{R}p}{\partial n} (\mathfrak{a}_p - \mathfrak{a}_w)& \nonumber \\ 
  - \frac{\partial \mathcal{R}i}{\partial n} (\mathfrak{a}_i - \mathfrak{a}_w) &= 0 
  \text{ in }\Omega \times (0,1],
  \\ 
   \mathfrak{a}_n(\mathbf{x}, 1) + (\obs p - p_{\text{d}}) \frac{\partial \obs p}{\partial n} 
  + (\obs n - n_{\text{d}}) \frac{\partial \obs n}{\partial n} 
  + (\obs l - l_{\text{d}}) \frac{\partial \obs l}{\partial n} & = 0 
   \text{ in }\Omega,
 \\ 
   - \partial_t \mathfrak{a}_o  + \frac{\partial \gamma}{\partial o} ( 1- n) (i \mathfrak{a}_i + p \mathfrak{a}_p - (i+p) \mathfrak{a}_n)
  + \frac{\partial \alpha}{\partial o} p (1 - i) (\mathfrak{a}_p - \mathfrak{a}_i) 
  &- \frac{\partial \beta}{\partial o} i ( 1 - p) (\mathfrak{a}_p - \mathfrak{a}_i) \nonumber \\
  + \delta_c p \mathfrak{a}_o + \delta_s w \mathfrak{a}_o
  - \frac{\partial \mathcal{R}p}{\partial o} (\mathfrak{a}_p - \mathfrak{a}_w) 
  - \frac{\partial \mathcal{R}i}{\partial o} (\mathfrak{a}_i - \mathfrak{a}_w) &= 0 
   \text{ in }\Omega \times (0,1],
 \\ 
   \mathfrak{a}_o(\mathbf{x}, 1) &= 0 \text{ in }\Omega ,
 \\ 
   - \partial_t \mathfrak{a}_w + \frac{\partial \mathcal{D} i}{\partial w} (\mathfrak{a}_w - \mathfrak{a}_i)
  - \delta_c (1 -o) \mathfrak{a}_o
  + \frac{\partial \mathcal{R}p}{\partial w} (\mathfrak{a}_w - \mathfrak{a}_p)& \nonumber \\ 
  + \frac{\partial \mathcal{R}i}{\partial w} (\mathfrak{a}_w - \mathfrak{a}_i) &= 0 
   \text{ in }\Omega \times (0,1],
  \\ 
   \mathfrak{a}_w(\mathbf{x}, 1) + (\obs p - p_{\text{d}}) \frac{\partial \obs p}{\partial w} 
  + (\obs n - n_{\text{d}}) \frac{\partial \obs n}{\partial w} 
  + (\obs l - l_{\text{d}}) \frac{\partial \obs l}{\partial w} &= 0 
   \text{ in }\Omega,
\end{align}
\end{subequations}
And finally, the inversion (gradient) equations for the parameters can be evaluated as, 
\begin{subequations}
  \label{eq:gradeq}
\begin{align}
  \frac{\partial \mathcal{L}}{\partial \kappa} & = \int_{0}^{1} \int_{\Omega} 
  - \frac{\partial \mathcal{D}i}{\partial \kappa} (\mathfrak{a}_i - \mathfrak{a}_w) ~\diff \Omega ~\diff t \\ 
  \frac{\partial \mathcal{L}}{\partial \rho} & = \int_{0}^{1} \int_{\Omega} 
  - \frac{\partial \mathcal{R}p}{\partial \rho}(\mathfrak{a}_p - \mathfrak{a}_w) 
  - \frac{\partial \mathcal{R}i}{\partial \rho}(\mathfrak{a}_p - \mathfrak{a}_w) ~\diff \Omega ~\diff t \\ 
  \frac{\partial \mathcal{L}}{\partial \beta_0} & = \int_{0}^{1} \int_{\Omega} 
  \frac{\partial \beta}{\partial \beta_0} i (1 - p) (\mathfrak{a}_i - \padj) ~\diff \Omega ~\diff t \\ 
  \frac{\partial \mathcal{L}}{\partial \alpha_0} & = \int_{0}^{1} \int_{\Omega} 
  \frac{\partial \alpha}{\partial \alpha_0} p (1 - i) (\padj - \iadj) ~\diff \Omega ~\diff t \\ 
  \frac{\partial \mathcal{L}}{\partial \gamma_0} & = \int_{0}^{1} \int_{\Omega} 
  \frac{\partial \gamma}{\gamma_0} (1-n)(p \padj - p \nadj + i \iadj - i \nadj) ~\diff \Omega ~\diff t \\
  \frac{\partial \mathcal{L}}{\partial \delta_c} & = \int_{0}^{1} \int_{\Omega} 
  \oadj o p ~\diff \Omega ~\diff t \\ 
  \frac{\partial \mathcal{L}}{\partial \delta_s} & = \int_{0}^{1} \int_{\Omega} 
  -\oadj (1 -o) w ~\diff \Omega ~\diff t \\ 
  \frac{\partial \mathcal{L}}{\partial \oinv} & = \int_{0}^{1} \int_{\Omega} 
  \frac{\partial \gamma}{\partial \ohyp} (1-n) (p \padj - p \nadj + i \iadj - i \nadj) ~\diff \Omega ~\diff t \\ 
  \frac{\partial \mathcal{L}}{\partial \ohyp} & = \int_{0}^{1} \int_{\Omega} 
  \frac{\partial \gamma}{\partial \ohyp} (1 - n) (p \padj - p \nadj + i \iadj - i \nadj) ~\diff \Omega ~\diff t \\ 
  \frac{\partial \mathcal{L}}{\partial \ith} & = \int_{\Omega} 
  -(\obs  l(\mathbf{x},1) - l_{\text{d}}) (1 - \obs p - \obs n) \hs'(i - \ith) ~\diff x 
\end{align}
\end{subequations}
where $\hs'$ is the derivative of the Heaviside function. 

To compute the gradient, it is necessary to first solve the forward problem, as outlined in \Cref{eq:forward_problem_1D}, then solve the adjoint equations, as per \Cref{eq:adjoint_equations_1D}, and then computing the gradient using \Cref{eq:gradeq}. 
To empirically study the ill-posedness of the 1D inversion problem, we compute the Hessian of the optimization setup through central differencing of the gradients computed in \Cref{eq:gradeq}. 

In \Cref{sssec:ill_posedness}, the results of the ill-posedness analysis are presented and the spectrum of eigenvalues of the Hessian for samples ($\mathbf{q}$) from coefficients space is illustrated in \Cref{fig:eig_spec}.
To overcome this challenge, a regularization method is proposed in \Cref{alg:reg} to penalize these ill-posed directions in the objective function. The method begins by uniformly sampling from the model coefficient space and estimating the Hessian at the minimum objective function value for each sample. The samples with smaller eigenvalues are deemed to contain less valuable information and, thus, their corresponding eigenvectors are weighted accordingly. A singular value decomposition is then performed to evaluate the weight of ill-posed directions among different coefficient combinations. The singular values are depicted in \Cref{fig:sing_values}. The right singular directions are identified and weighted based on the inverse of the singular values ($\Uill$) to mitigate directions with small singular values in the objective function. These directions with smaller singular values are more ill-posed and, therefore, require more penalization in the objective function. 
\begin{algorithm}[h]
  \caption{Algorithm to compute the regularization}\label{alg:reg}
  \begin{algorithmic}[1]    
    \Procedure{}{$n$, $\mathbf{b}_l$, $\mathbf{b}_u$} \Comment{Get the number of samples $n$, lower bound vector $\mathbf{b}_l$ and upper bounds vector $\mathbf{b}_u$}
    \State $\{\mathbf{q}_i\}_{i=1}^{n} \sim U(\mathbf{b}_l, \mathbf{b}_u)$ \Comment{Take uniform samples of parameters}
    \State $\mathbf{D} \leftarrow \mathbf{0}$
    \For{$i=1 \dots n$}
      \State $\mathbf{H} \leftarrow \text{Hess}(\mathbf{q}_i)$ \Comment{Compute the Hessian}
      \State $\Lambda, \mathbf{Q} \leftarrow \text{EvalDecomp}(\mathbf{H})$  \Comment{Compute Evals. and Evecs. of Hessian}
      \State $\mathbf{D} \leftarrow \mathbf{D} \cup \Lambda^{\frac{1}{2}} \mathbf{Q} $ \Comment{Store the weighted ill-posed directions}
    \EndFor
    \State $\mathbf{\Sigma}, \mathbf{U} \leftarrow \text{SVD}(\mathbf{D})$ \Comment{SVD on ill-posed directions}
    \State $\mathbf{\Uill} \leftarrow \mathbf{U} \mathbf{\Sigma}^{-1}$ \Comment{Weighting the ill-posed directions}
    \EndProcedure
  \end{algorithmic}
\end{algorithm}
\subsection{3D Inversion}
As mentioned, the inverse problem requires three parameters to characterize the tumor progression: 
\begin{itemize*}
  \item Healthy brain anatomy $\Omega$ to describe $w(\mathbf{x},0), g(\mathbf{x},0), f(\mathbf{x},0)$
  \item Initial condition of proliferative tumor cells $p_0$
  \item Model coefficients vector $\mathbf{q}$
\end{itemize*}. 
In the following, we explain how we estimate each of these parameters. 
\subsubsection{Brain anatomy}
\label{sssec:brain_anatomy}
The process of deriving the healthy brain anatomy ($\Omega$) involves multiple steps. First, an MRI scan is affine registered to a template atlas. Then, a neural network \cite{lipkova2019personalized} trained on the BraTS dataset \cite{bakas2018identifying} is used to segment the tumor into proliferative, necrotic, and edema regions. Next, ten normal brain scans with known healthy segmentation are registered to the patient's scan using \cite{avants2009advanced}. An ensemble-based approach is then used to combine the atlas-based segmentations for the different healthy tissue labels \cite{gholami2019novel,subramanian2022ensemble}. We combine the healthy tissue labels and tumor segmentation as a single brain segmentation of the patient. Finally, the tumorous regions are replaced with white matter to estimate the healthy brain anatomy $\Omega$.
\subsubsection{IC Inversion}
\label{sssec:IC}
To estimate IC of proliferative tumor cells, we use the single-species IC inversion proposed in \cite{scheufele2020fully}. This algorithm uses an adjoint based method for reaction-diffusion PDE model constrained with sparsity and $\max$ constraints. Our current inverse IC scheme deviates from the \cite{scheufele2020fully} algorithm by limiting the support region to the necrotic areas only, as our simulations in \cite{subramanian2019simulation} revealed that the IC always present in the necrotic regions. 
\subsubsection{Coefficients Inversion}
To invert for the model coefficients $\mathbf{q}$ in 3D, we use the regularization term computed from the 1D model as sampling the Hessian is very expensive in 3D. 
As the structure of the nonlinear operators is the same, our hypothesis is the 1D regularization operator can be used to stabilize the 3D reconstruction. One contribution of this paper is to verify the hypothesis using numerical experiments with synthetically constructed data. 
Given the IC ($p_0$) and brain anatomy ($\Omega$) estimated from \Cref{sssec:brain_anatomy,sssec:IC}, we can solve the following optimization problem with constraints:
\begin{align}
  \min_{\mathbf{q}} \mathcal{J} = \frac{1}{2} \| \obs p(\mathbf{x},1) - p_{\text{d}} \|^2_{L_2(\Omega)} & + \frac{1}{2} \| \obs n(\mathbf{x},1) - n_{\text{d}} \|^2_{L_2(\Omega)} + \frac{1}{2} \| \obs l(\mathbf{x},1) - l_{\text{d}} \|^2_{L_2(\Omega)} + \frac{\lambda}{2} \| \Uill^T \mathbf{q} \|^2 \\ 
  \text{subject to } & \mathcal{F}(p_0,\mathbf{q}) \text{ in } \Omega \times (0, 1] 
\end{align}

\section{Results}
\label{sec:results}
In this study, we evaluate our algorithm on synthetic cases and clinical dataset. We carry out the 3D inversions using the Frontera and Lonestar system at the Texas Advanced Computing Center (TACC) at The University of Texas at Austin. Our solver is written in both C++ and Python. The forward solver uses C++ with PETSc, AccFFT, PnetCDF. The inverse problem is solved using the covariance matrix adaptation evolution strategy (CMA-ES) in Python \cite{hansen2019pycma, hansen1996adapting}. For 1D test-cases, we use Python for the forward problem with the same CMA-ES optimizer.

\subsubsection{Numerical parameters}
We list all numerical parameters used in our solver in \Cref{tab:numerical_parameters}. We now explain a justification for our choice:
\begin{enumerate}[label=(\roman*)]
 \item For the multi-species forward solver, we use the same setting and numerics of \cite{subramanian2019simulation}.  
 \item For IC inversion, we use the solver from our previous studies \cite{scheufele2020fully,subramanian2020did}. The single-species reaction diffusion model with $\ell_0$ constraint is used to solve a single-snapshot inverse problem iteratively. The IC is estimated within necrotic regions and the sparsity parameter is fixed at $5$ to allow up to 5 voxels to be considered for IC. The experiments on choosing this value can be found in \cite{scheufele2020fully}. 
 \item We determine the optimal value for the regularization parameter $\lambda$ by evaluating the quality of data reconstruction with $20\%$ additive noise. We select this noise percentage because previous studies have demonstrated that the tumor segmentation tool \cite{bakas2018identifying} achieves an accuracy of approximately $80\%$. 
 \item For the inversion of model coefficients $\mathbf{q}$, we employ CMA-ES with a relative function tolerance of $1\times 10^{-8}$ and an absolute function tolerance of $1\times 10^{-3}$. Our analysis indicates that reducing these tolerance values further does not result in improved outcomes. 
 \item The population size of the CMA-ES is $16$ as a balance between evaluation costs and having a diverse population for covariance sampling. 
 \item We set the initial guess for the CMA-ES solver to the midpoint of the parameter bounds and the initial variance to half the parameter range. Despite experimenting with higher variances and different initial guesses, we did not observe any improvement in solver performance.
 \item We use a shape factor $\omega$ of $32$ for all the Heaviside functions except for $\hs(i-\ith)$ to achieve a smooth transition. For the edema region, a smooth Heaviside function in $\hs(i-\ith)$ can lead to high spatial diffusivity of infiltrative cells $i$. To prevent this, we use a larger shape factor ($\omega=256$) for the edema region compared to the others.
\end{enumerate}

 \subsubsection{Performance measures}
 Our main goal of the numerical experiments is to show how the proposed regularization can stabilize the multi-species inversion and then use the same regularization to perform a full clinical inversion of the algorithm. First, we show how the regularization helps the inversion in synthetic cases and then combine it with IC inversion on clinical scans. In this regard, we define the following metrics to quantify our reconstructions, 
 \begin{enumerate}[label=(\roman*)]
  \item relative error in model coefficient $\iota$:
    \begin{align}
      e_{\iota} = \frac{|\iota^{\rec}-\iota^{*}|}{| \iota^{*}|},
    \end{align}
     where $\iota^{*}$ is ground truth coefficient and $\iota^{\rec}$ is the reconstructed value. We denote the model coefficients vector with $\mathbf{q}$ 
     and $e_{\mathbf{q}}$ denotes the average relative error of model coefficients. 
   \item relative error in final species reconstruction $s(t=1)$ where $s$ can be $n$, $p$ or $i$:
   \begin{align}
    \mu_{s, L_{2}} = \frac{\| s^{\rec}(1) - s^{*}(1) \|_{L_2(\Omega)}}{\|s^{*}(1)\|_{L_2(\Omega)}},
   \end{align} 
   where $s^{\rec}(1)$ is the final reconstructed species, $s^{*}(1)$ is the final species concentration grown with ground truth model coefficients $\mathbf{q}^{*}$. 
   \item relative error in initial tumor concentration:
   \begin{align}
     \mu_{0, L_2} = \frac{\| p^{\rec}_0 - p^{*}_0 \|_{L_2(\Omega)}}{\|p^{*}_0\|_{L_2(\Omega)}},
   \end{align}
   where $p^{\rec}(0)$ is the reconstructed IC and $p^{*}(0)$ is the ground truth IC. 
   \item relative $\ell_1$ error in initial tumor parameterization ($p_0$):
   \begin{align}
      \mu_{0, \ell_1} = \frac{\| p_0^{\rec} - p_0^{*} \|_{\ell_1}}{\| p_0^{*}\|_{\ell_1}},
   \end{align}
   where $p_0^{\rec}$ is the reconstructed tumor IC and $p_0^{*}$ is the ground truth IC. 
   \item relative $\ell_1$ reconstruction for projected IC along $a$ axis:
   \begin{align}
    \mu_{0, \ell_1}^{a} =
    \frac{\| \int_{a} (p_0^{\rec} - p_0^{*}) \diff a \|_{\ell_1}}{\| \int_{a} p_0^{*} \diff a \|_{\ell_1}}
   \end{align}
   where $a$ represents the axis which can be $x,y$ or $z$. 
   \item reconstruction dice score for an observed species $s$:
   \begin{align}
      \pi_{s} = \frac{2 | s^{\rec} \bigcap s^{*} | }{ | s^{\rec} | + | s^{*} |}
   \end{align}
   where $| s^{\rec} |$ and $| s^{*} |$ are the cardinality number of the observed species from inversion and data, respectively and $| s^{\rec} \bigcap s^{*} |$ is the cardinality number of the intersection between inverted species $s^{\rec}$  and species from data $s^{*}$. 
   \item reconstruction dice score for tumor core ($\obs p \cup \obs n$) denoted as $\text{tc}$
   \begin{align}
    \pi_{tc} = \frac{2 | \text{tc}^{\rec} \bigcap \text{tc}^{*} | }{ | \text{tc}^{\rec} | + | \text{tc}^{*} |}
   \end{align}
   where $| \text{tc}^{\rec} |$ and $| \text{tc}^{*} |$ are the cardinality number of tumor core from inversion and data, respectively and $| \text{tc}^{\rec} \bigcap \text{tc}^{*} |$ is the cardinality number of the intersection between inverted tumor core $\text{tc}^{\rec}$  and tumor core from data $ \text{tc}^{*}$. Note that we denote the dice score from single-species model with $\pi_{tc}^{\text{ss}}$ and dice score from multi-species model with $\pi_{tc}^{\text{ms}}$.       
 \end{enumerate}
 \begin{table}[!ht]
  \centering
  \caption{Numerical parameters for tumor inversion and forward solver. Initial guess is mid-range for model coefficients. We refer to \cite{scheufele2020fully} for IC inversion parameters in more details.}
  \resizebox{0.7\textwidth}{!}{
      \begin{tabular}{cc}
          \toprule 
          Parameter & Value \\ 
          \midrule
          Spatial discretization for 1D inversion & $256$ \\
          Spatial discretization for 3D inversion & $160^3$ \\
          Time-step of the forward solvers in 1D and 3D   & $\num{2e-2}$ \\ 
          Shape factor $\omega$ in $\obs c$, $\obs n, \obs p, \alpha, \beta$ and $\gamma$  & $ 32$ \\ 
          Shape factor $\omega$ in $\hs(i-\ith)$ & $256$ \\ 
          Epsilon used in finite differencing of Hessian, $\epsilon$ & $\num{1e-6}$ \\ 
          Regularization parameter in 1D inversion, $\lambda$ & $\num{1e-4}$ \\ 
          Regularization parameter in 3D inversion, $\lambda$ & $\num{1e-5}$ \\ 
          Relative function tolerance & $\num{1e-7}$ \\ 
          Absolute function tolerance & $\num{1e-3}$ \\ 
          Population size of the CMA-ES solver & $16$ \\ 
          \bottomrule
        \end{tabular}
        }        
        \label{tab:numerical_parameters}
\end{table}  

 \subsection{1D test-case}
 In this section, first, we discuss the ill-posedness of the problem in 1D and then, we present the inversion results. We present the analysis of the 1D inversion assuming that the IC is known and we only need to invert for the model coefficients.  
 \subsubsection{Ill-posedness analysis}
 \label{sssec:ill_posedness}
\begin{figure}
  \centering
  \begin{subfigure}{0.45\textwidth}
    \centering
    \includegraphics[width=\textwidth]{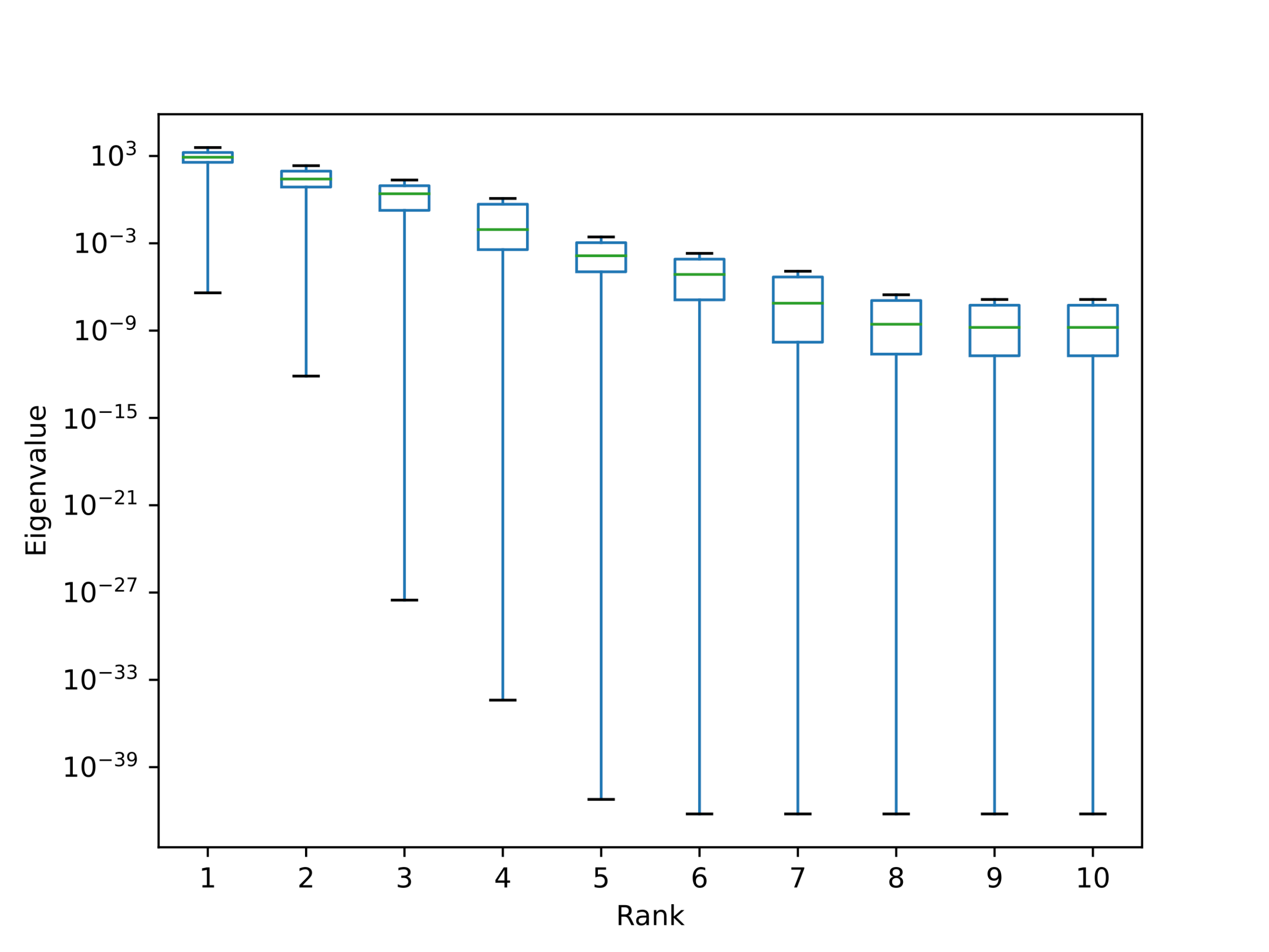}
    \caption{Spectrum of eigenvalues for Hessian}
    \label{fig:eig_spec}
  \end{subfigure}
  \hfill
  \begin{subfigure}{0.45\textwidth}
    \centering
    \includegraphics[width=\textwidth]{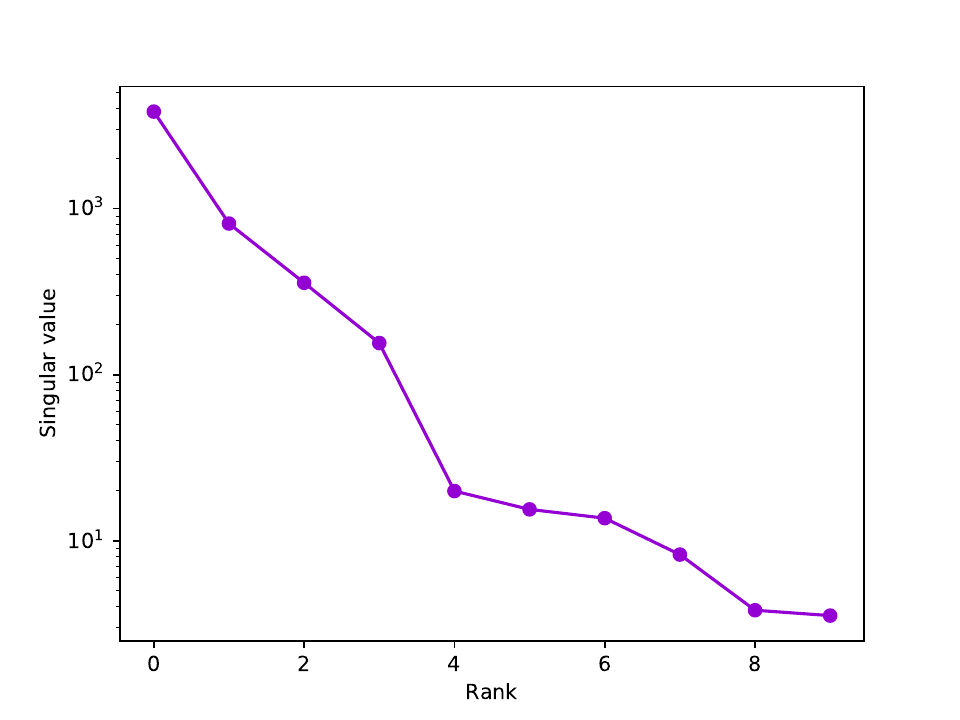}
    \caption{Singular values of weighted eigenvectors}    
    \label{fig:sing_values}
  \end{subfigure}
  \vfill
  \begin{subfigure}{0.5\textwidth}
    \vspace{4mm}
    \centering
    \includegraphics[width=\textwidth]{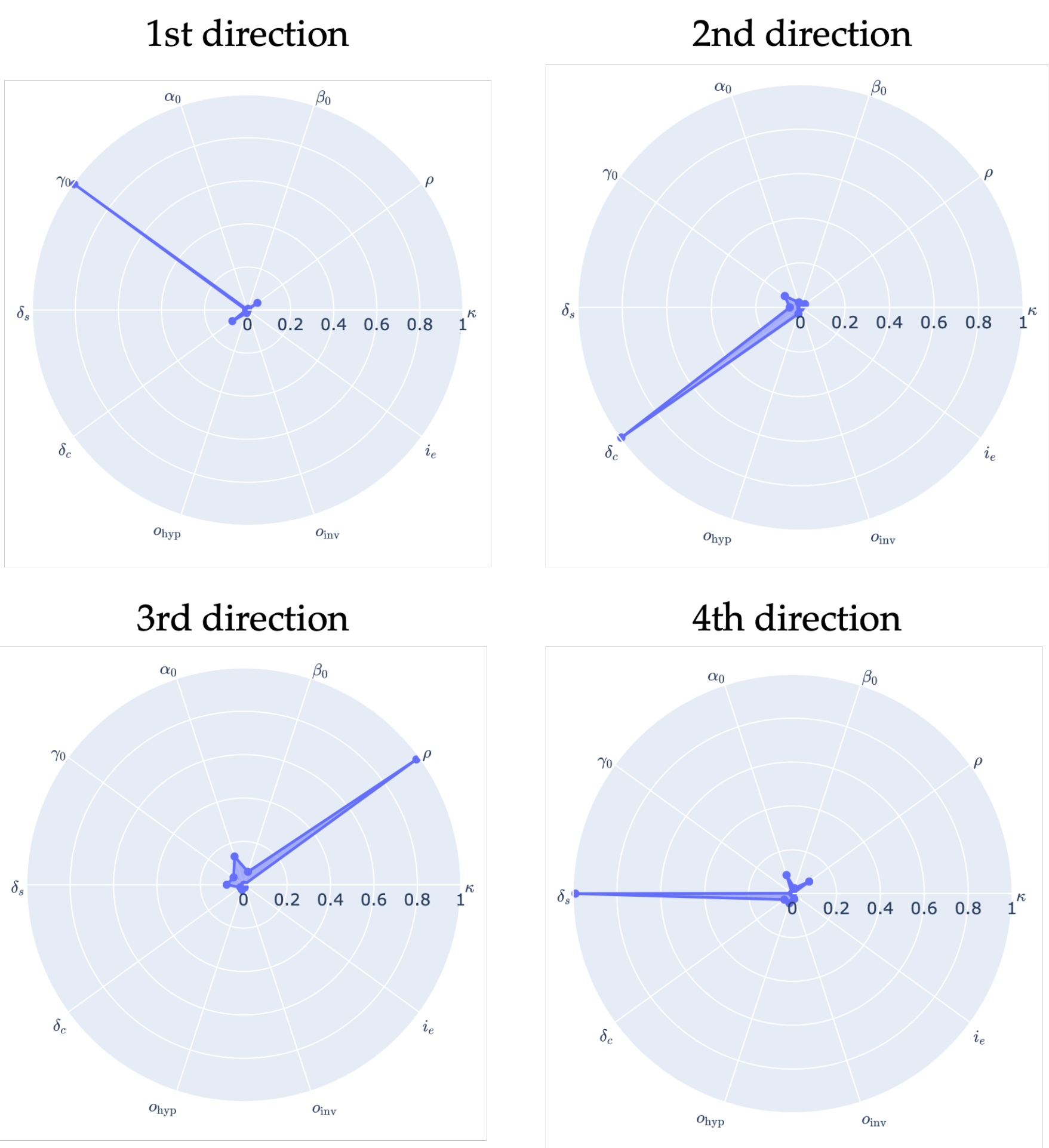}
    \caption{Four most ill-posed directions computed using \Cref{alg:reg}}    
    \label{fig:directions}
  \end{subfigure}
  \vspace{4mm}
  \caption{Ill-posedness analysis of multi-species model in 1D. The spectrum of Hessian eigenvalues for 3500 parameter combinations is presented in Figure~\ref{fig:eig_spec}, with each eigenvalue normalized such that the max eigenvalue is 1. The singular values calculated in Algorithm~\ref{alg:reg} are shown in Figure~\ref{fig:sing_values}. We depict the absolute value of last four directions of the computed regularization term ($\Uill$) in \Cref{fig:directions} corresponding to the four smallest singular values. The exponential decrease of eigenvalues in Figure~\ref{fig:eig_spec} indicates high ill-posedness of the problem. Additionally, the wide range of small eigenvalues suggests varying degrees of ill-posedness among samples, which is accounted for in Algorithm~\ref{alg:reg}. Figure~\ref{fig:sing_values} demonstrates how different directions impact the ill-posedness, with smaller singular values indicating greater ill-posedness. Thus, using the inverse of the singular values can lead to improved inversion results. The results in \Cref{fig:directions} clearly indicates ill-posedness for high values $\gamma_0$ and a combination of high values of $\rho$ and low values of $\delta_c$ contributes to ill-posedness.}
\end{figure}

 To address the ill-posedness of the problem, we evaluate the Hessian for 3500 samples of $\mathbf{q}$ uniformly derived from coefficients range and perform an eigenvalue decomposition. The eigenvalue spectrum, displayed in \Cref{fig:eig_spec}, exhibits an exponential decrease, indicating that the problem is highly ill-posed. The wide range of eigenvalues suggests that different combinations of coefficients contribute to the problem's ill-posedness with varying degrees of difficulty to reconstruct the coefficients. 

 Following \Cref{alg:reg}, we further analyze the ill-posed directions by performing singular value decomposition (SVD) on weighted eigenvectors. The result of the SVD is displayed in \Cref{fig:sing_values} through the plot of the singular values. The smaller singular value suggests higher ill-posedness of the corresponding coefficient combinations.
 In \Cref{fig:directions}, we have illustrated the last four directions of the computed $\Uill$, which correspond to the 4 smallest singular values from \Cref{fig:sing_values}. It is crucial to consider the weight for penalizing the $\gamma_0$ parameter in the first direction. This penalization corresponds to the necrotic species and accounts for the fact that if a location is assigned to the necrotic region in both the data and observation, increasing the $\gamma_0$ coefficient will not reduce the error as the objective function is insensitive to it. Furthermore, we have observed through our simulations that a combination of high values of $\rho$ and low values of $\delta_c$ in the second direction contribute to the ill-posedness.
 
 \subsubsection{Inversion analysis in 1D with noise} 
 In this section, we present results from inverting 1D synthetic cases, where exact IC is known and focus is on finding model coefficients $\mathbf{q}$. To test inversion stability, we generate species concentrations ($s^*$) using ground truth model coefficients $\mathbf{q}^{*}$, and we generate the noisy species $\hat{s}^* = s^{*} (1 + \textit{noise})$ where \textit{noise} has a zero-mean gaussian distribution with variance $\sigma$ ($\mathcal{N}(0, \sigma)$). We then apply the observation operators to generate the noisy observed data. We vary $\sigma$ to generate different average noise levels among all species $p,n$ and $i$. This noisy data is then input for the inversion algorithm and tested with and without regularization to evaluate stability and performance of reconstruction.

\begin{figure}
  \centering
  \includegraphics[width=\textwidth]{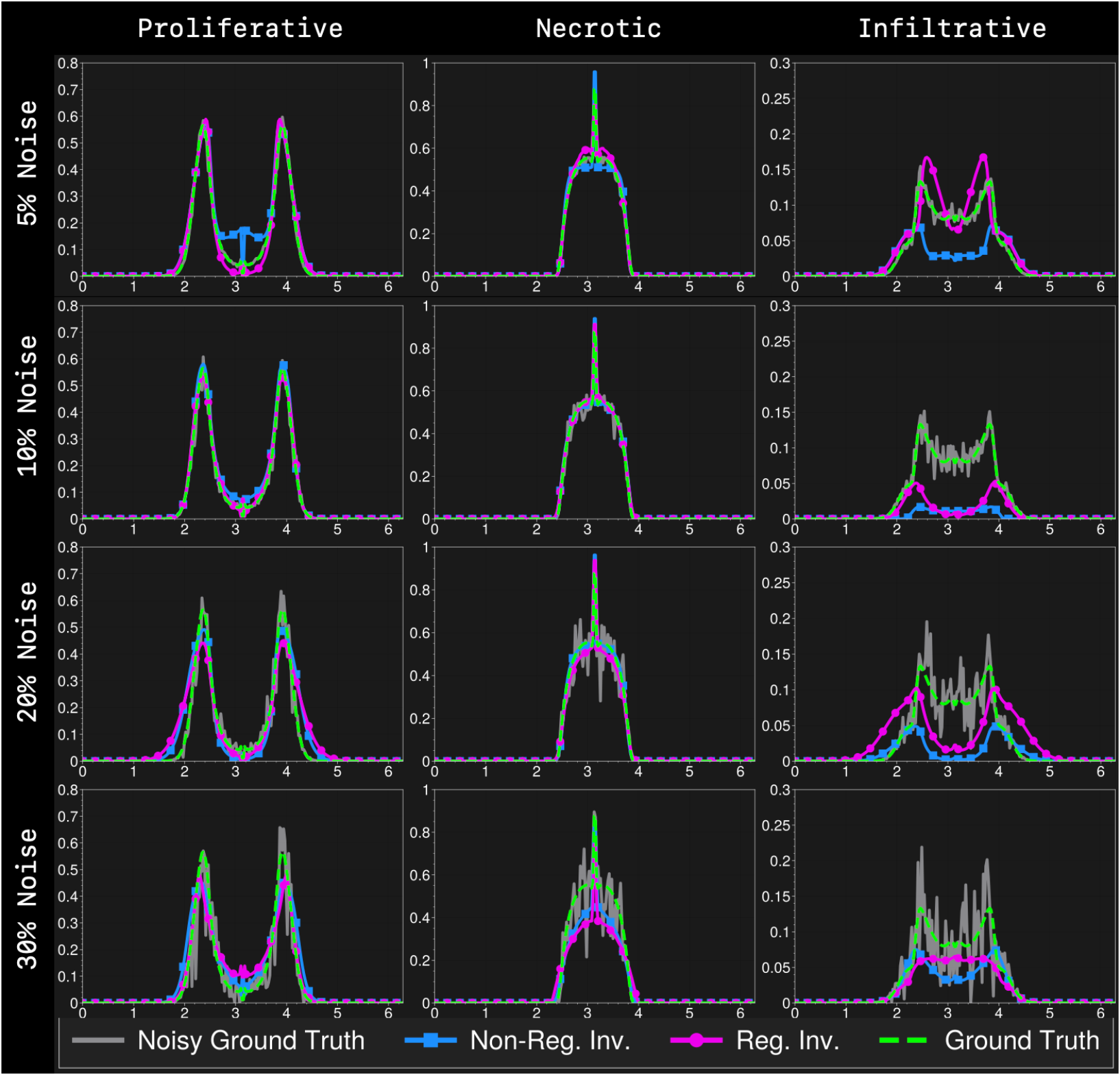}
  \caption{1D inversion results of synthetic data with additive Gaussian noise ($5\%,10\%,20\%,30\%$). Here we show the reconstructed species concentration at $T=1$. First to fourth rows show relative noise levels. Plots in first to third columns show proliferative, necrotic, and infiltrative species, respectively. The green, gray, magenta, and blue lines represent the underlying ground truth species (before observation denoted as "Ground Truth"), ground truth species with noise ("Noisy Ground Truth"), reconstructed species without regularization ("Non-Reg. Inv."), and reconstructed species with regularization ("Reg. Inv."), respectively. Regularized inversion is observed to better match the underlying species even when the data is perturbed by noise, particularly in the case of infiltrative cells where deviations are greatly reduced with regularization.}
  \label{fig:1D_noise}
 \end{figure}

 In \Cref{fig:1D_noise}, we show reconstructed species concentrations (proliferative $p$, infiltrative $i$, necrotic $n$) for synthetic cases with different levels of independent Gaussian noise added to each signal. The results with and without regularization are displayed. For a quantitative comparison, see \Cref{tab:1D_noise_err} where regularization generally improves reconstruction performance compared to non-regularized results. 
  In \Cref{tab:1D_noise_err}, we also observe that regularization avoids parameter inversion problems, as seen in the case with $10\%$ noise where some parameters diverge. \Cref{tab:1D_noise_err_params} reports inverted coefficients for each parameter with and without regularization. 
  Overall, the inversion scheme with regularization is more stable in terms of species reconstruction and parameter inversion. 

 Further, we examine the influence of varying coefficient combinations on the performance of the inversion process. To achieve this objective, we conduct experiments on $150$ synthetic data sets contaminated with $20\%$ Gaussian noise and the results are presented in \Cref{tab:1D_noise_err_params_samples}. The reconstruction measures are calculated both for regularized and non-regularized inversion cases. Our results indicate that the use of regularization leads to a significant reduction in the average relative error of the coefficients ($e_{\mathbf{q}}$), with an improvement of approximately $30\%$. This improvement is particularly noticeable for the coefficients $\gamma_0, \beta_0,$ and $\alpha_0$.

\begin{table}[!ht]
  \centering
  \caption{Evaluation of 1D inversion scheme for the model parameters with synthetic data assuming known initial condition for the proliferative tumor cells. Gaussian noise is added to species to generate noisy data. Dice scores of observed species $s$ ($\pi_{s}$) and relative $L_2$ error of true species with inverted species $s$ ($\mu_{s, L_2}$) are reported. Average of relative errors of model coefficients $\mathbf{q}$ ($e_{\mathbf{p}}$) is also reported. Regularization improves reconstruction error, especially for species and model coefficients, though not for all noise levels. Nevertheless, regularization stabilizes the inversion scheme with respect to non-regularized inversion where the deviation of errors is significant. The dice coefficients in the unregularized case are better due to overfitting to the data}
  \resizebox{\textwidth}{!}{
      \begin{tabular}{cc|ccc|ccc|c}
          \toprule 
          Case & Noise (\%) & $\pi_{n}$ & $\pi_{p}$ & $\pi_{l}$ &  $\mu_{n, L_{2}}$ & $\mu_{p, L_{2}}$ & $\mu_{i, L_{2}}$ & $e_{\mathbf{q}}$                                 
          \\
          \midrule
          Non-Reg & $5$ & ${1.65e-02}$ & ${9.36e-01}$ & ${4.07e-02}$ & ${7.74e-02}$ & ${2.43e-01}$ & ${5.99e-01}$ & ${6.64e-01}$ \\ 
          %\rowcolor{gray!30}
          Reg & $5$ & ${9.86e-01}$ & ${9.36e-01}$ & ${5.80e-01}$ & ${8.01e-02}$ & ${1.45e-01}$ & ${3.69e-01}$ & ${5.23e-01}$ \\           
          \midrule
          Non-Reg & $10$ & ${1.77e-02}$ & ${9.34e-01}$ & ${9.89e-02}$ & ${4.91e-02}$ & ${1.11e-01}$ & ${8.73e-01}$ & ${1.36e+00}$ \\ 
          %\rowcolor{gray!30}
          Reg & $10$ & ${9.83e-01}$ & ${9.29e-01}$ & ${4.26e-01}$ & ${5.79e-02}$ & ${8.12e-02}$ & ${7.59e-01}$ & ${6.86e-01}$ \\           
          \midrule
          Non-Reg & $20$ & ${1.50e-02}$ & ${8.92e-01}$ & ${5.59e-02}$ & ${6.04e-02}$ & ${2.89e-01}$ & ${8.07e-01}$ & ${6.76e-01}$ \\ 
          %\rowcolor{gray!30}
          Reg & $20$ & ${9.70e-01}$ & ${8.99e-01}$ & ${5.13e-01}$ & ${9.86e-02}$ & ${3.25e-01}$ & ${6.87e-01}$ & ${7.98e-01}$ \\           
          \midrule
          Non-Reg & $30$ & ${6.85e-02}$ & ${8.22e-01}$ & ${4.43e-02}$ & ${2.41e-01}$ & ${2.64e-01}$ & ${4.96e-01}$ & ${4.65e-01}$ \\
          %\rowcolor{gray!30}
          Reg & $30$ & ${9.17e-01}$ & ${8.32e-01}$ & ${3.71e-01}$ & ${3.17e-01}$ & ${2.56e-01}$ & ${4.14e-01}$ & ${4.74e-01}$ \\           
          \toprule
      \end{tabular}
  }
  \label{tab:1D_noise_err}
\end{table}          
 \subsection{3D test-cases}
 In this section, we analyze the 3D inversion. Results for cases with known brain anatomy ($\Omega$) and IC ($p_0$) is presented first, where our goal is to estimate growth coefficients $\mathbf{q}$ only. The next case involves inverting for both IC and growth coefficients. Finally, results for real clinical data from BraTS20 training dataset \cite{bakas2018identifying} are presented.
 \subsubsection{Artificial tumor with known IC and anatomy}
 \label{sssec:AT_KnownIC}
\begin{figure}
  \centering
  \includegraphics[width=\textwidth]{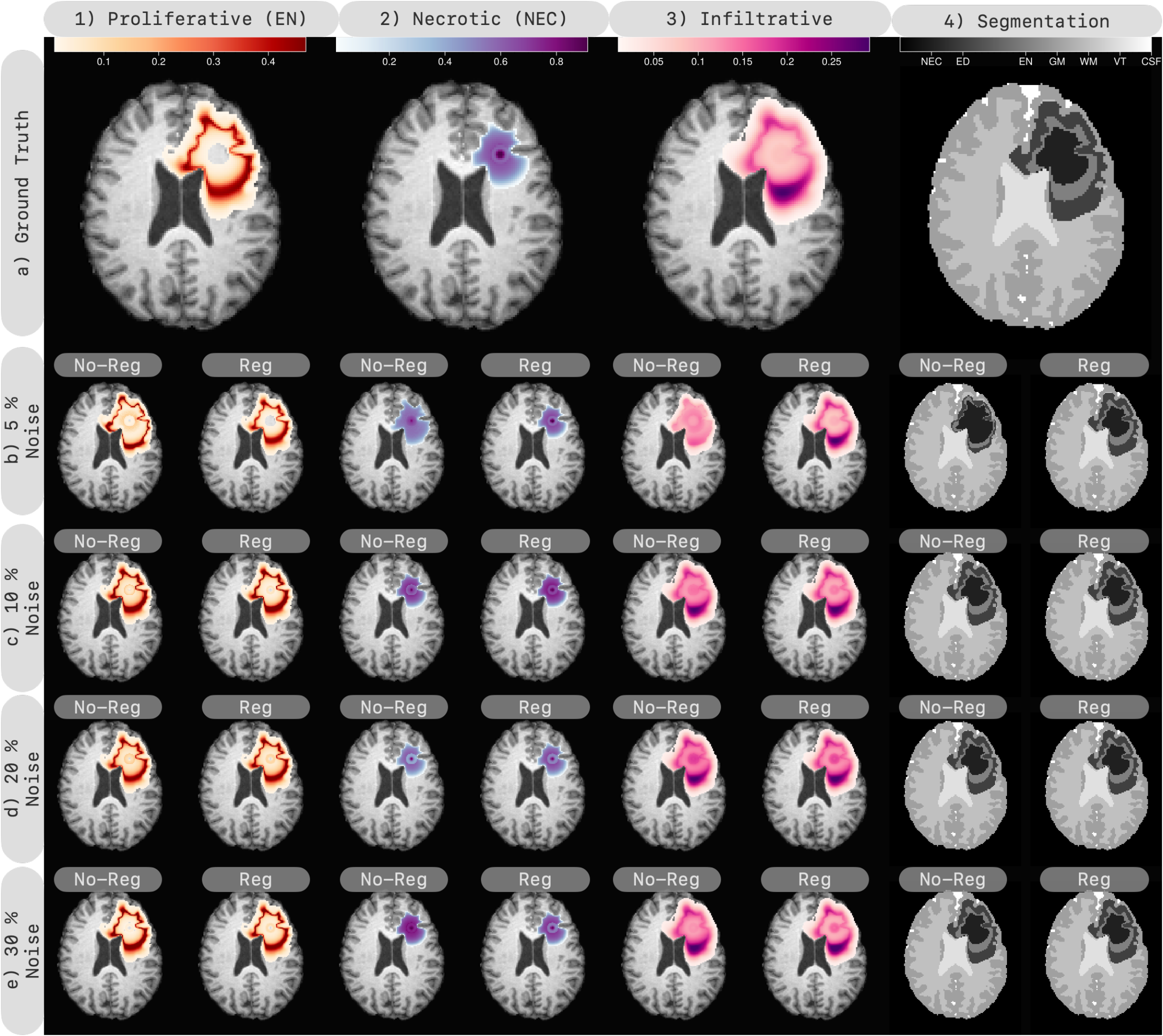}
  \caption{Evaluating 3D inversion using synthetic data with known IC and brain anatomy. We start by generating a tumor test-case and showing its species and segmentation in the first row (a). We then add $5\%$, $10\%$, $20\%$ and $30\%$  gaussian noise to the species and apply an observation operator to get noisy data. We perform inversion with and without regularization on this data. The inversion results are displayed in rows 2-4, with the enhancing, necrotic, and infiltrative concentrations in columns 1-3, and the segmentation in the last column. The regularized inversion yields acceptable segmentation for all noise levels, whereas the non-regularized case shows  shows more variability for different noise levels.}
  \label{fig:3d_noise}
\end{figure}

 In this section, we evaluate inversion for growth parameters only. 
 Given a ground truth $\mathbf{q}^*$, we solve the forward problem to generate data. We test different noise levels of $5\%, 10\%, 20\%$, and $30\%$ and present the reconstruction with and without regularization in \Cref{fig:3d_noise}. We report the results in \Cref{tab:3D_noise_err}.

 As seen in \Cref{tab:3D_noise_err} and \Cref{fig:3d_noise}, regularization significantly improves the reconstruction error $e_{\mathbf{q}}$ in $5\%$ noise, 
 while improving the reconstruction at higher noise levels. Although dice scores and segmentation do not show significant improvement, the underlying species reconstruction, not directly observed from data, is improved.

 \begin{table}[!ht]
  \centering
  \caption{Numerical evaluation of 3D inversion using synthetic data with known anatomy and IC. We compare the dice score of observed species $\pi_{s}$, the relative $L_2$ error of species $\mu_{s, L_2}$, and average relative error of coefficients $e_{\mathbf{q}}$ for noise levels of $5\%$, $10\%$, $20\%$, and $30\%$ with and without regularization. Regularization improves performance, especially for the $5\%$ noise level in terms of coefficient error, although dice scores do not vary drastically. The significance can be seen in the relative $L_2$ error of underlying species.}
  \resizebox{\textwidth}{!}{
      \begin{tabular}{cc|ccc|ccc|c}
          \toprule 
          Case & Noise (\%) & $\pi_{n}$ & $\pi_{p}$ & $\pi_{l}$ &  $\mu_{n, L_{2}}$ & $\mu_{p, L_{2}}$ & $\mu_{i, L_{2}}$ & $e_{\mathbf{q}}$\\
          \midrule
          Non-Reg & 5& $\num{1.98e-01}$ & $\num{5.72e-01}$ & $\num{5.87e-01}$ & $\num{9.95e-01}$ & $\num{9.36e-01}$ & $\num{4.21e-01}$ & $\num{1.30e+00}$ \\ 
          %\rowcolor{gray!30}
          Reg & 5& $\num{9.70e-01}$ & $\num{9.81e-01}$ & $\num{9.84e-01}$ & $\num{4.41e-02}$ & $\num{5.66e-02}$ & $\num{2.91e-02}$ & $\num{2.57e-01}$ \\ 
          Non-Reg & 10& $\num{9.23e-01}$ & $\num{9.53e-01}$ & $\num{9.52e-01}$ & $\num{2.14e-01}$ & $\num{2.51e-01}$ & $\num{1.07e-01}$ & $\num{7.58e-01}$ \\ 
          %\rowcolor{gray!30}
          Reg & 10& $\num{9.35e-01}$ & $\num{9.58e-01}$ & $\num{9.61e-01}$ & $\num{2.07e-01}$ & $\num{2.83e-01}$ & $\num{8.90e-02}$ & $\num{6.69e-01}$ \\ 
          Non-Reg & 20& $\num{9.57e-01}$ & $\num{9.72e-01}$ & $\num{9.73e-01}$ & $\num{1.01e-01}$ & $\num{1.44e-01}$ & $\num{1.39e-01}$ & $\num{9.42e-01}$ \\ 
          %\rowcolor{gray!30}
          Reg & 20& $\num{9.49e-01}$ & $\num{9.67e-01}$ & $\num{9.67e-01}$ & $\num{1.18e-01}$ & $\num{1.36e-01}$ & $\num{1.21e-01}$ & $\num{4.81e-01}$ \\ 
          Non-Reg & 30& $\num{9.24e-01}$ & $\num{9.53e-01}$ & $\num{9.44e-01}$ & $\num{1.94e-01}$ & $\num{3.40e-01}$ & $\num{1.42e-01}$ & $\num{3.76e-01}$ \\ 
          %\rowcolor{gray!30}
          Reg & 30& $\num{9.54e-01}$ & $\num{9.71e-01}$ & $\num{9.70e-01}$ & $\num{1.06e-01}$ & $\num{1.83e-01}$ & $\num{1.13e-01}$ & $\num{3.30e-01}$ \\ 
          \toprule
      \end{tabular}
  }
  \label{tab:3D_noise_err}
\end{table}     
\subsubsection{Artificial tumor with unknown IC and anatomy}
\label{sssec:AT_UknownIC}
In this section, we present the results for two-stage inversion process for both IC ($p_0$) and model coefficients ($\mathbf{q}$). Our solution involves first the tumorous regions are filled with white matter and the healthy brain is estimated and we invert for the IC using a single-species model, as detailed in \cite{subramanian2020did,scheufele2020fully}. Once the IC is estimated, we estimate the growth model coefficients. 

\begin{figure}
  \centering
  \includegraphics[width=\textwidth]{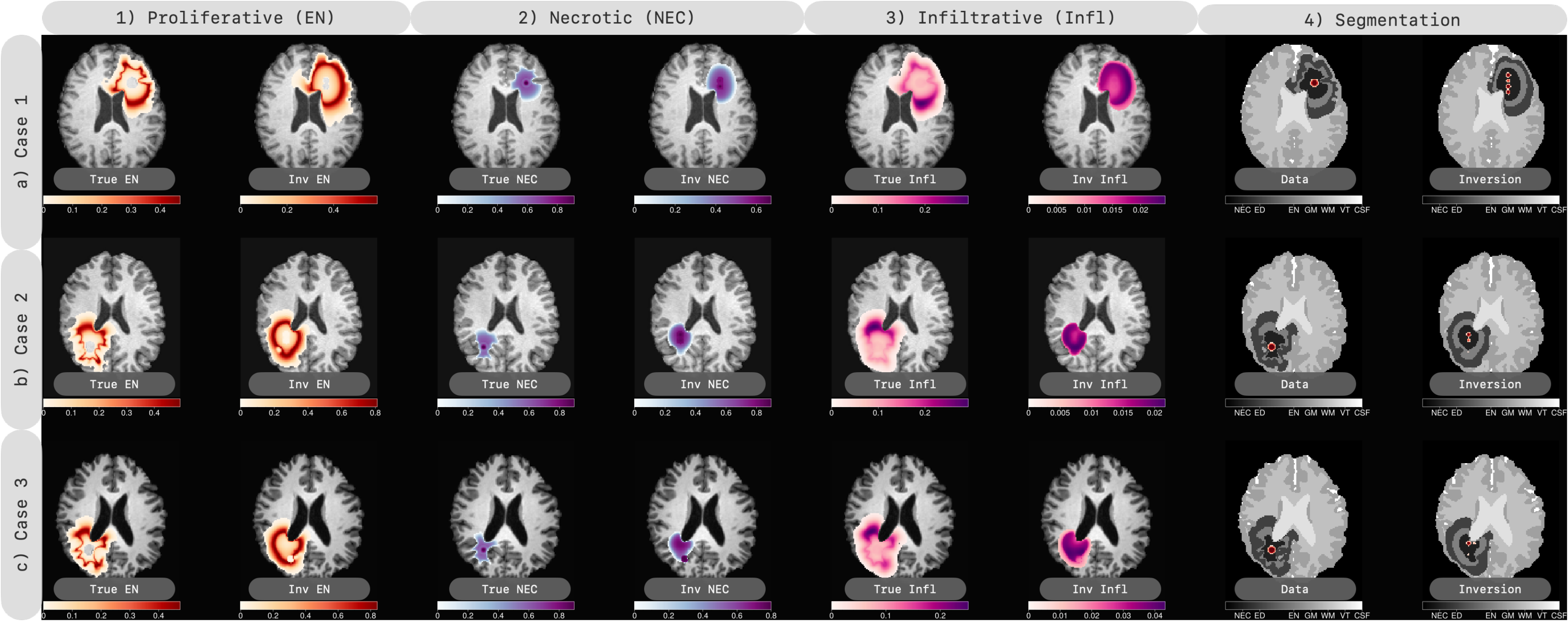}
  \caption{3D inversion result with synthetic test-cases generated with an unknown IC and brain anatomy. Different normal brains with no mass deformation are used to generate the synthetic cases. In each row, we present a single test-case, with tumor species (Proliferative (EN), necrotic (NEC), and infiltrative cells) from the synthetic cases and inversion displayed in the first, second, and third columns, respectively. The given data and segmentation computed by the inversion algorithm are depicted in the fourth column. Moreover, we depicted the projected ground truth IC and estimated IC in the segmentation for data and inversion, respectively, Although the IC is an estimate, the reconstruction from the inversion is good as evidenced by the qualitative observation. The challenge of reconstructing multi-species arises due to the unknown initial brain anatomy and isotropic migration of tumor in all directions due to assigned white matter label.}
  \label{fig:IC_syn}
\end{figure}

Three cases are tested with generated data and different parameter combinations. The inversion results are depicted in \Cref{fig:IC_syn} and the corresponding numerical measures are reported in \Cref{tab:IC_syn_err}. The individual model coefficients and their relative error are also included in the \Cref{tab:IC_syn_err}.

As we can observe, we demonstrate that even though the IC estimation is derived from another model, it is still capable of reconstructing the observed tumor species from the data. Our test cases, as shown in \Cref{fig:IC_syn}, indicate that the IC estimation provides a relatively good estimate for the initial condition. The estimation of the IC for a reaction-diffusion tumor growth model is exponentially ill-posed, making it a challenging task, as noted in prior studies such as \cite{subramanian2020did,scheufele2020fully}.

As shown in \Cref{tab:IC_syn_err}, our solver not only achieves a reconstruction comparable to that of the single-species model but also provides additional information on the underlying species. Our method enables the characterization of observed species and the estimation of non-observable ones from the given data.
As demonstrated in \Cref{sssec:AT_KnownIC}, our algorithm exhibits good convergence once a suitable estimate of IC is obtained. Notably, the performance of our method is influenced by the brain tumor anatomy, with the diffusion patterns significantly impacted by the tissue types in the test cases. Our model simplifies the problem by estimating the tumor growth model, given the challenge of accurately determining brain anatomy due to mass deformation. Furthermore, as illustrated in \Cref{fig:IC_syn}, the diffusion of proliferative cells is primarily in all directions due to the allocation of white matter in the brain anatomy, while the data is influenced by the gray matter.

\begin{table}[!ht]
  \centering
  \caption{Quantitative evaluation of 3D inversion with unknown IC using synthetic data. The accuracy of species reconstruction is measured by the dice score $\pi_s$, and the relative $L_2$ reconstruction of each species is reported as $\mu_{s, L_2}$. The relative $\ell_1$ and $L_2$ errors for the IC are indicated as $\mu_{0,L_2}$ and $\mu_{0, \ell_1}$, respectively, and the relative error in model coefficients is denoted by $e_{\mathbf{q}}$. We also report the projected error of IC estimation denoted as $\mu_{0, \ell_1}^{x}$, $\mu_{0, \ell_1}^{y}$ and $\mu_{0, \ell_1}^{z}$ in $yz$, $xz$ and $xy$ planes, respectively. The results suggest that the dice scores exhibit good accuracy in reconstructing the observed species. However, the combination of IC errors and coefficient errors significantly affects the relative error in the parameters, pointing to the need for further improvement.}
  \resizebox{\textwidth}{!}{
      \begin{tabular}{c|ccccc|ccc|ccccc|c}
          \toprule 
          Case & $\pi_{\text{tc}}^{\text{ms}}$ & $\pi_{\text{tc}}^{\text{ss}}$ & $\pi_{p}$ & $\pi_{n}$ & $\pi_{l}$ &  $\mu_{n, L_{2}}$ & $\mu_{p, L_{2}}$ & $\mu_{i, L_{2}}$  & $\mu_{0, L_2}$ & $\mu_{0, \ell_1}$ & $\mu_{0, \ell_1}^{x}$ &$\mu_{0, \ell_1}^{y}$ & $\mu_{0, \ell_1}^{z}$ & $e_{\mathbf{q}}$\\
          \midrule
          Case 1 & $\num{8.43e-01}$ & $\num{7.73e-01}$ & $\num{6.46e-01}$ & $\num{7.59e-01}$ & $\num{7.85e-01}$ & $\num{9.78e-01}$ & $\num{5.08e-01}$ & $\num{9.00e-01}$ & $\num{1.18e+00}$ & $\num{1.60e+00}$ & $\num{1.14e+00}$ & $\num{1.12e+00}$ & $\num{1.05e+00}$ &  $\num{4.28e+00}$\\
          Case 2 & $\num{7.95e-01}$ & $\num{7.35e-01}$ & $\num{5.43e-01}$ & $\num{7.74e-01}$ & $\num{5.74e-01}$ & $\num{1.28e+00}$ & $\num{8.78e-01}$ & $\num{9.48e-01}$ & $\num{1.10e+00}$ & $\num{1.33e+00}$ & $\num{1.06e+00}$ & $\num{6.82e-01}$ & $\num{1.06e+00}$ & $\num{1.36e+00}$\\
          Case 3 & $\num{8.19e-01}$ & $\num{6.78e-01}$ & $\num{6.13e-01}$ & $\num{7.93e-01}$ & $\num{7.20e-01}$ & $\num{1.26e+00}$ & $\num{7.64e-01}$ & $\num{8.94e-01}$ & $\num{1.09e+00}$ & $\num{1.27e+00}$ & $\num{1.06e+00}$ & $\num{1.04e+00}$ & $\num{1.07e+00}$ & $\num{8.28e-01}$ \\
          \toprule
      \end{tabular}
  }
  \label{tab:IC_syn_err}
\end{table}

\subsection{Inversion on two clinical data}
\begin{figure}
  \centering
  \includegraphics[width=\textwidth]{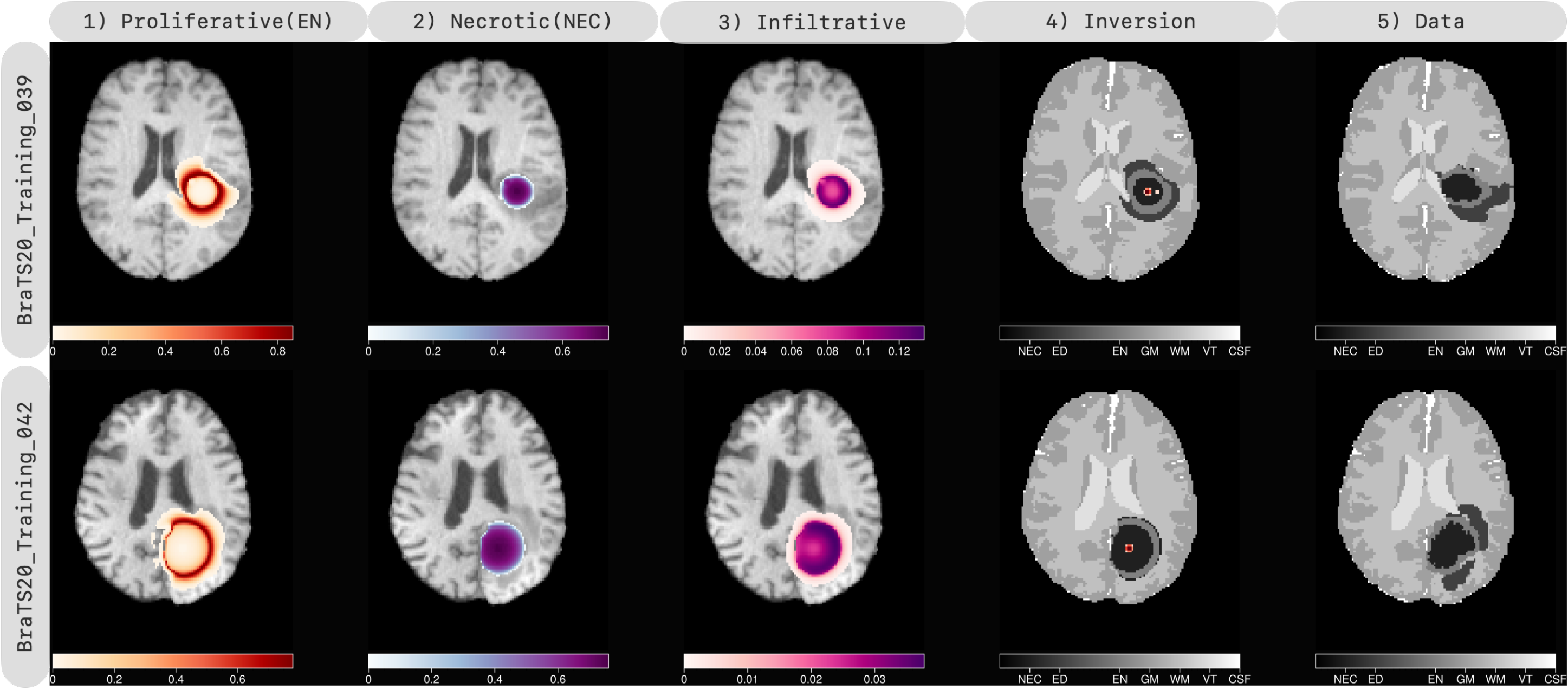}
  \caption{We conduct a multi-species inversion of the BraTS20\_Training clinical dataset, selecting two patients, using the MRI scan segmentation as input. Our results are presented in a tabular format, with each row showing the inversion and data for each patient. The first three columns depict the inverted species (enhancing (EN), necrotic (NEC), and infiltrative cells), the fourth column shows the observed species from the inversion, and the last column displays the data. Our results demonstrate that, despite only having a single snapshot from each patient, we can qualitatively reconstruct the species distribution with reasonable accuracy.}
  \label{fig:real_IC}
\end{figure}

In this section, we demonstrate the application of our inversion scheme to two patients from the BraTS20 dataset \cite{bakas2018identifying}. The procedure for reconstruction is outlined in \Cref{sssec:AT_UknownIC} and involves first reconstructing IC of the tumor, followed by estimation of the growth parameters. The reconstructed species are visualized in \Cref{fig:real_IC}, with the projected IC estimate shown in the inverted segmentation. Numerical measures of the results can be found in \Cref{tab:IC_clinical} and \Cref{tab:3D_clinical_params}. It is important to note that the only information available from the patients is the segmentation of the observed species from MRI scans, and direct comparison of the underlying species concentrations is not possible.

Our solver achieves a reconstruction accuracy comparable to that of the single-species model, similar to the findings in \Cref{sssec:AT_UknownIC}. Furthermore, we are able to reconstruct species data and quantify infiltrative tumor cells. However, the reconstruction may be affected by IC estimation using the single-species model, particularly for BraTS20\_Training\_039, as evidenced by the presence of a considerable amount of necrotic tissue in the initial condition location. Furthermore, the anisotropic progression of edema in BraTS20\_Training\_039 presents a challenging task for the solver to identify infiltrative cells. A similar trend can be observed for BraTS20\_Training\_042, where the edema does not uniformly diffuse. Our edema model is a simplified representation, used in  \cite{subramanian2019simulation, saut2014multilayer}. Additionally, the results may be influenced by the species' volume in the data, as our unweighted objective function prioritizes reducing errors for the dominant species, i.e. necrotic tissue in BraTS20\_Training\_042. 

\begin{table}[!ht]
  \centering
  \caption{The results of the multi-species model inversion using clinical data from the BraTS20 dataset are presented in the table. The performance of the model is measured by the dice score, which serves as a quantitative measure of the accuracy of the reconstructed tumor species. The dice score for the tumor core reconstruction by the multi-species model is denoted as $\pi^{\text{ms}}{\text{tc}}$, while the dice score for the tumor core reconstruction by the single-species model is denoted as $\pi^{\text{ss}}{\text{tc}}$. Additionally, the dice scores for the reconstruction of necrotic, proliferative, and edema regions are denoted as $\pi_{n}$, $\pi_{p}$, and $\pi_{l}$, respectively. As is evident from the results, our solver exhibits similar performance in reconstructing the tumor core compared to the single-species model. Additionally, it provides the reconstruction of the individual species, offering added value in comparison to the single-species approach.}
  \resizebox{0.7\textwidth}{!}{
      \begin{tabular}{c|ccccc}
          \toprule 
          Patient's ID & $\pi^{\text{ms}}_{\text{tc}}$ & $\pi^{\text{ss}}_{\text{tc}}$ & $\pi_{n}$ & $\pi_{p}$ & $\pi_{l}$ \\
          \midrule
          BraTS20\_Training\_039 & $\num{8.23e-01}$ & $\num{8.37e-01}$ & $\num{6.94e-01}$ & $\num{6.18e-01}$ & $\num{4.24e-01}$ \\
          BraTS20\_Training\_042 & $\num{9.01e-01}$ & $\num{8.92e-01}$ & $\num{7.25e-01}$ & $\num{5.41e-01}$ & $\num{0.00e+00}$ \\
          \toprule
      \end{tabular}
  }
  \label{tab:IC_clinical}
\end{table}

\section{Conclusion}
\label{sec:conclusions}
In this study, we present an inverse solver for tumor growth model from one time snapshot  and a regularization scheme to allow stable reconstruction of model parameters. 
. The solver estimates the initial condition of the multi-species model and growth coefficients. The input data used is the segmentation of the tumor obtained from MRI scans, and we design a single-snapshot inverse solver for the multi-species model.

Our results demonstrate the efficacy of this regularization in improving the estimation of the growth coefficients in one-dimensional and three-dimensional settings, even in the presence of high relative noise. Moreover, we combine the regularization terms with our prior work on tumor initial condition estimation and apply the inversion algorithm to both synthetic and clinical cases, revealing the reconstruction's stability for the underlying species concentration. Our solver not only estimates the tumor core in a range comparable to that of the single-species model, but also provides additional information on the reconstruction of other species. Therefore, we can estimate the non-observable species for the tumor growth model.

The solver is implemented using a scalable sampling optimization scheme, demonstrating a reasonable computational time. Although the estimation of the initial condition could be improved to enhance the model coefficients and reconstructions, our results provide valuable insights into the inverse problem for multi-species tumor growth models.

In the future, we plan to conduct a comprehensive evaluation of the inversion scheme on a large number of clinical images to observe the correlation between the model coefficients and clinical data outcomes, such as survival rate. Given the complexity of the biophysical model, there are numerous opportunities for future extensions, including the inclusion of mass deformation for more efficient computation of the inverse problem.

\section{Acknowledgements}

\bibliography{refs.bib}

\begin{thebibliography}{10}
\providecommand \doibase [0]{http://dx.doi.org/}%

\bibitem{subramanian2019simulation}
Subramanian S, Gholami A, Biros G. Simulation of glioblastoma growth using a 3D
  multispecies tumor model with mass effect. {\it Journal of mathematical
  biology.} 2019\string;79(3)\string:941--967.

\bibitem{subramanian2020did}
Subramanian S, Scheufele K, Mehl M, Biros G. Where did the tumor start? An
  inverse solver with sparse localization for tumor growth models. {\it Inverse
  problems.} 2020\string;36(4)\string:045006.

\bibitem{scheufele2020fully}
Scheufele K, Subramanian S, Biros G. Fully automatic calibration of
  tumor-growth models using a single mpMRI scan. {\it IEEE transactions on
  medical imaging.} 2020\string;40(1)\string:193--204.

\bibitem{subramanian2020multiatlas}
Subramanian S, Scheufele K, Himthani N, Biros G. Multiatlas calibration of
  biophysical brain tumor growth models with mass effect. In: Springer.
  2020\string:551--560.

\bibitem{subramanian2022ensemble}
Subramanian S, Ghafouri A, Scheufele K, Himthani N, Davatzikos C, Biros G.
  Ensemble inversion for brain tumor growth models with mass effect. {\it IEEE
  Transactions on Medical Imaging.} 2022.

\bibitem{hogea2007modeling}
Hogea C, Davatzikos C, Biros G. Modeling glioma growth and mass effect in 3D MR
  images of the brain. In: Springer.  2007\string:642--650.

\bibitem{hogea2008brain}
Hogea C, Davatzikos C, Biros G. Brain--Tumor interaction biophysical models for
  medical image registration. {\it SIAM Journal on Scientific Computing.}
  2008\string;30(6)\string:3050--3072.

\bibitem{jbabdi2005simulation}
Jbabdi S, Mandonnet E, Duffau H, et al. Simulation of anisotropic growth of
  low-grade gliomas using diffusion tensor imaging. {\it Magnetic Resonance in
  Medicine: An Official Journal of the International Society for Magnetic
  Resonance in Medicine.} 2005\string;54(3)\string:616--624.

\bibitem{mang2012biophysical}
Mang A, Toma A, Schuetz TA, et al. Biophysical modeling of brain tumor
  progression: From unconditionally stable explicit time integration to an
  inverse problem with parabolic PDE constraints for model calibration. {\it
  Medical Physics.} 2012\string;39(7Part1)\string:4444--4459.

\bibitem{hormuth2018mechanically}
Hormuth DA, Eldridge SL, Weis JA, Miga MI, Yankeelov TE. Mechanically coupled
  reaction-diffusion model to predict glioma growth: methodological details.
  {\it Cancer Systems Biology: Methods and Protocols.} 2018\string:225--241.

\bibitem{saut2014multilayer}
Saut O, Lagaert JB, Colin T, Fathallah-Shaykh HM. A multilayer grow-or-go model
  for GBM: effects of invasive cells and anti-angiogenesis on growth. {\it
  Bulletin of mathematical biology.} 2014\string;76(9)\string:2306--2333.

\bibitem{swanson2011quantifying}
Swanson KR, Rockne RC, Claridge J, Chaplain MA, Alvord~Jr EC, Anderson AR.
  Quantifying the role of angiogenesis in malignant progression of gliomas: in
  silico modeling integrates imaging and histology. {\it Cancer research.}
  2011\string;71(24)\string:7366--7375.

\bibitem{gholami2019novel}
Gholami A, Subramanian S, Shenoy V, et al. A novel domain adaptation framework
  for medical image segmentation. In: Springer.  2019\string:289--298.

\bibitem{scheufele2020image}
Scheufele K, Subramanian S, Mang A, Biros G, Mehl M. Image-driven biophysical
  tumor growth model calibration. {\it SIAM journal on scientific computing: a
  publication of the Society for Industrial and Applied Mathematics.}
  2020\string;42(3)\string:B549.

\bibitem{fathi2018characterization}
Fathi~Kazerooni A, Nabil M, Zeinali~Zadeh M, et al. Characterization of active
  and infiltrative tumorous subregions from normal tissue in brain gliomas
  using multiparametric MRI. {\it Journal of Magnetic Resonance Imaging.}
  2018\string;48(4)\string:938--950.

\bibitem{peeken2019deep}
Peeken JC, Molina-Romero M, Diehl C, et al. Deep learning derived tumor
  infiltration maps for personalized target definition in Glioblastoma
  radiotherapy. {\it Radiotherapy and Oncology.}
  2019\string;138\string:166--172.

\bibitem{gholami2016inverse}
Gholami A, Mang A, Biros G. An inverse problem formulation for parameter
  estimation of a reaction--diffusion model of low grade gliomas. {\it Journal
  of mathematical biology.} 2016\string;72(1)\string:409--433.

\bibitem{lipkova2019personalized}
Lipkov{\'a} J, Angelikopoulos P, Wu S, et al. Personalized radiotherapy design
  for glioblastoma: Integrating mathematical tumor models, multimodal scans,
  and bayesian inference. {\it IEEE transactions on medical imaging.}
  2019\string;38(8)\string:1875--1884.

\bibitem{hormuth2021image}
Hormuth DA, Al~Feghali KA, Elliott AM, Yankeelov TE, Chung C. Image-based
  personalization of computational models for predicting response of high-grade
  glioma to chemoradiation. {\it Scientific reports.}
  2021\string;11(1)\string:1--14.

\bibitem{ghafouri20233d}
Ghafouri A, Biros G. A 3D Inverse Solver for a Multi-species PDE Model of
  Glioblastoma Growth. In: Springer.  2023\string:51--60.

\bibitem{swanson2000quantitative}
Swanson KR, Alvord~Jr EC, Murray J. A quantitative model for differential
  motility of gliomas in grey and white matter. {\it Cell proliferation.}
  2000\string;33(5)\string:317--329.

\bibitem{gooya2012glistr}
Gooya A, Pohl KM, Bilello M, et al. GLISTR: glioma image segmentation and
  registration. {\it IEEE transactions on medical imaging.}
  2012\string;31(10)\string:1941--1954.

\bibitem{bakas2018identifying}
Bakas S, Reyes M, Jakab A, et al. Identifying the best machine learning
  algorithms for brain tumor segmentation, progression assessment, and overall
  survival prediction in the BRATS challenge. {\it arXiv preprint
  arXiv:1811.02629.} 2018.

\bibitem{ozisik2018inverse}
Ozisik MN. {\it Inverse heat transfer: fundamentals and applications}.
\newblock Routledge, 2018.

\bibitem{cheng2008quasi}
Cheng J, Liu J. A quasi Tikhonov regularization for a two-dimensional backward
  heat problem by a fundamental solution. {\it Inverse Problems.}
  2008\string;24(6)\string:065012.

\bibitem{zheng2014recovering}
Zheng GH, Wei T. Recovering the source and initial value simultaneously in a
  parabolic equation. {\it Inverse Problems.} 2014\string;30(6)\string:065013.

\bibitem{avants2009advanced}
Avants BB, Tustison N, Song G, others . Advanced normalization tools (ANTS).
  {\it Insight j.} 2009\string;2(365)\string:1--35.

\bibitem{hansen1996adapting}
Hansen N, Ostermeier A. Adapting arbitrary normal mutation distributions in
  evolution strategies: The covariance matrix adaptation. In: IEEE.
  1996\string:312--317.

\end{thebibliography}


@article{subramanian2019simulation,
  title={Simulation of glioblastoma growth using a 3D multispecies tumor model with mass effect},
  author={Subramanian, Shashank and Gholami, Amir and Biros, George},
  journal={Journal of mathematical biology},
  volume={79},
  number={3},
  pages={941--967},
  year={2019},
  publisher={Springer}
}
@article{fathi2018characterization,
  title={Characterization of active and infiltrative tumorous subregions from normal tissue in brain gliomas using multiparametric MRI},
  author={Fathi Kazerooni, Anahita and Nabil, Mahnaz and Zeinali Zadeh, Mehdi and Firouznia, Kavous and Azmoudeh-Ardalan, Farid and Frangi, Alejandro F and Davatzikos, Christos and Saligheh Rad, Hamidreza},
  journal={Journal of Magnetic Resonance Imaging},
  volume={48},
  number={4},
  pages={938--950},
  year={2018},
  publisher={Wiley Online Library}
}
@article{peeken2019deep,
  title={Deep learning derived tumor infiltration maps for personalized target definition in Glioblastoma radiotherapy},
  author={Peeken, Jan C and Molina-Romero, Miguel and Diehl, Christian and Menze, Bjoern H and Straube, Christoph and Meyer, Bernhard and Zimmer, Claus and Wiestler, Benedikt and Combs, Stephanie E},
  journal={Radiotherapy and Oncology},
  volume={138},
  pages={166--172},
  year={2019},
  publisher={Elsevier}
}
@article{gholami2016inverse,
  title={An inverse problem formulation for parameter estimation of a reaction--diffusion model of low grade gliomas},
  author={Gholami, Amir and Mang, Andreas and Biros, George},
  journal={Journal of mathematical biology},
  volume={72},
  number={1},
  pages={409--433},
  year={2016},
  publisher={Springer}
}
@article{saut2014multilayer,
  title={A multilayer grow-or-go model for GBM: effects of invasive cells and anti-angiogenesis on growth},
  author={Saut, Olivier and Lagaert, Jean-Baptiste and Colin, Thierry and Fathallah-Shaykh, Hassan M},
  journal={Bulletin of mathematical biology},
  volume={76},
  number={9},
  pages={2306--2333},
  year={2014},
  publisher={Springer}
}
@article{subramanian2020did,
  title={Where did the tumor start? An inverse solver with sparse localization for tumor growth models},
  author={Subramanian, Shashank and Scheufele, Klaudius and Mehl, Miriam and Biros, George},
  journal={Inverse problems},
  volume={36},
  number={4},
  pages={045006},
  year={2020},
  publisher={IOP Publishing}
}
@article{scheufele2020fully,
  title={Fully automatic calibration of tumor-growth models using a single mpMRI scan},
  author={Scheufele, Klaudius and Subramanian, Shashank and Biros, George},
  journal={IEEE transactions on medical imaging},
  volume={40},
  number={1},
  pages={193--204},
  year={2020},
  publisher={IEEE}
}
@inproceedings{subramanian2020multiatlas,
  title={Multiatlas calibration of biophysical brain tumor growth models with mass effect},
  author={Subramanian, Shashank and Scheufele, Klaudius and Himthani, Naveen and Biros, George},
  booktitle={International Conference on Medical Image Computing and Computer-Assisted Intervention},
  pages={551--560},
  year={2020},
  organization={Springer}
}
@article{mang2020integrated,
  title={Integrated biophysical modeling and image analysis: application to neuro-oncology},
  author={Mang, Andreas and Bakas, Spyridon and Subramanian, Shashank and Davatzikos, Christos and Biros, George},
  journal={Annual review of biomedical engineering},
  volume={22},
  pages={309},
  year={2020},
  publisher={NIH Public Access}
}
@article{hormuth2021image,
  title={Image-based personalization of computational models for predicting response of high-grade glioma to chemoradiation},
  author={Hormuth, David A and Al Feghali, Karine A and Elliott, Andrew M and Yankeelov, Thomas E and Chung, Caroline},
  journal={Scientific reports},
  volume={11},
  number={1},
  pages={1--14},
  year={2021},
  publisher={Nature Publishing Group}
}
@article{hormuth2019calibrating,
  title={Calibrating a predictive model of tumor growth and angiogenesis with quantitative MRI},
  author={Hormuth, David A and Jarrett, Angela M and Feng, Xinzeng and Yankeelov, Thomas E},
  journal={Annals of biomedical engineering},
  volume={47},
  number={7},
  pages={1539--1551},
  year={2019},
  publisher={Springer}
}
@article{lipkova2019personalized,
  title={Personalized radiotherapy design for glioblastoma: Integrating mathematical tumor models, multimodal scans, and bayesian inference},
  author={Lipkov{\'a}, Jana and Angelikopoulos, Panagiotis and Wu, Stephen and Alberts, Esther and Wiestler, Benedikt and Diehl, Christian and Preibisch, Christine and Pyka, Thomas and Combs, Stephanie E and Hadjidoukas, Panagiotis and others},
  journal={IEEE transactions on medical imaging},
  volume={38},
  number={8},
  pages={1875--1884},
  year={2019},
  publisher={IEEE}
}
@article{subramanian2022ensemble,
  title={Ensemble inversion for brain tumor growth models with mass effect},
  author={Subramanian, Shashank and Ghafouri, Ali and Scheufele, Klaudius and Himthani, Naveen and Davatzikos, Christos and Biros, George},
  journal={IEEE Transactions on Medical Imaging},
  year={2022},
  publisher={IEEE}
}
@article{bakas2018identifying,
  title={Identifying the best machine learning algorithms for brain tumor segmentation, progression assessment, and overall survival prediction in the BRATS challenge},
  author={Bakas, Spyridon and Reyes, Mauricio and Jakab, Andras and Bauer, Stefan and Rempfler, Markus and Crimi, Alessandro and Shinohara, Russell Takeshi and Berger, Christoph and Ha, Sung Min and Rozycki, Martin and others},
  journal={arXiv preprint arXiv:1811.02629},
  year={2018}
}
@article{prusiner1990transgenetic,
  title={Transgenetic studies implicate interactions between homologous PrP isoforms in scrapie prion replication},
  author={Prusiner, Stanley B and Scott, Michael and Foster, Dallas and Pan, Keh-Ming and Groth, Darlene and Mirenda, Carol and Torchia, Marilyn and Yang, Shu-Lian and Serban, Dan and Carlson, George A and others},
  journal={Cell},
  volume={63},
  number={4},
  pages={673--686},
  year={1990},
  publisher={Elsevier}
}
@article{hansen2020cma,
  title={CMA-ES/pycma: r3. 0.3},
  author={Hansen, Nikolaus and Akimoto, Y and Baudis, P},
  journal={Cma-es/pycma: r3. 0.3},
  year={2020}
}
@inproceedings{hansen1996adapting,
  title={Adapting arbitrary normal mutation distributions in evolution strategies: The covariance matrix adaptation},
  author={Hansen, Nikolaus and Ostermeier, Andreas},
  booktitle={Proceedings of IEEE international conference on evolutionary computation},
  pages={312--317},
  year={1996},
  organization={IEEE}
}
@article{himthani2022claire,
  title={CLAIRE—Parallelized Diffeomorphic Image Registration for Large-Scale Biomedical Imaging Applications},
  author={Himthani, Naveen and Brunn, Malte and Kim, Jae-Youn and Schulte, Miriam and Mang, Andreas and Biros, George},
  journal={Journal of Imaging},
  volume={8},
  number={9},
  pages={251},
  year={2022},
  publisher={MDPI}
}
@inproceedings{himthani2023claire,
  title={CLAIRE: Scalable multi-GPU algorithms for diffeomorphic image registration},
  author={Himthani, Naveen Nares and Brunn, Malte and Schulte, Miriam and Biros, George and Mang, Andreas},
  booktitle={2023 Joint Mathematics Meetings (JMM 2023)},
  year={2023},
  organization={AMS}
}
@inproceedings{hogea2007modeling,
  title={Modeling glioma growth and mass effect in 3D MR images of the brain},
  author={Hogea, Cosmina and Davatzikos, Christos and Biros, George},
  booktitle={Medical Image Computing and Computer-Assisted Intervention--MICCAI 2007: 10th International Conference, Brisbane, Australia, October 29-November 2, 2007, Proceedings, Part I 10},
  pages={642--650},
  year={2007},
  organization={Springer}
}
@article{hogea2008brain,
  title={Brain--Tumor interaction biophysical models for medical image registration},
  author={Hogea, Cosmina and Davatzikos, Christos and Biros, George},
  journal={SIAM Journal on Scientific Computing},
  volume={30},
  number={6},
  pages={3050--3072},
  year={2008},
  publisher={SIAM}
}
@article{jbabdi2005simulation,
  title={Simulation of anisotropic growth of low-grade gliomas using diffusion tensor imaging},
  author={Jbabdi, Sa{\^a}d and Mandonnet, Emmanuel and Duffau, Hugues and Capelle, Laurent and Swanson, Kristin Rae and P{\'e}l{\'e}grini-Issac, M{\'e}lanie and Guillevin, R{\'e}my and Benali, Habib},
  journal={Magnetic Resonance in Medicine: An Official Journal of the International Society for Magnetic Resonance in Medicine},
  volume={54},
  number={3},
  pages={616--624},
  year={2005},
  publisher={Wiley Online Library}
}
@article{mang2012biophysical,
  title={Biophysical modeling of brain tumor progression: From unconditionally stable explicit time integration to an inverse problem with parabolic PDE constraints for model calibration},
  author={Mang, Andreas and Toma, Alina and Schuetz, Tina A and Becker, Stefan and Eckey, Thomas and Mohr, Christian and Petersen, Dirk and Buzug, Thorsten M},
  journal={Medical Physics},
  volume={39},
  number={7Part1},
  pages={4444--4459},
  year={2012},
  publisher={Wiley Online Library}
}
@article{hormuth2018mechanically,
  title={Mechanically coupled reaction-diffusion model to predict glioma growth: methodological details},
  author={Hormuth, David A and Eldridge, Stephanie L and Weis, Jared A and Miga, Michael I and Yankeelov, Thomas E},
  journal={Cancer Systems Biology: Methods and Protocols},
  pages={225--241},
  year={2018},
  publisher={Springer}
}
@article{pham2012density,
  title={Density-dependent quiescence in glioma invasion: instability in a simple reaction--diffusion model for the migration/proliferation dichotomy},
  author={Pham, Kara and Chauviere, Arnaud and Hatzikirou, Haralambos and Li, Xiangrong and Byrne, Helen M and Cristini, Vittorio and Lowengrub, John},
  journal={Journal of biological dynamics},
  volume={6},
  number={sup1},
  pages={54--71},
  year={2012},
  publisher={Taylor \& Francis}
}
@article{swanson2011quantifying,
  title={Quantifying the role of angiogenesis in malignant progression of gliomas: in silico modeling integrates imaging and histology},
  author={Swanson, Kristin R and Rockne, Russell C and Claridge, Jonathan and Chaplain, Mark A and Alvord Jr, Ellsworth C and Anderson, Alexander RA},
  journal={Cancer research},
  volume={71},
  number={24},
  pages={7366--7375},
  year={2011},
  publisher={AACR}
}
@article{scheufele2020image,
  title={Image-driven biophysical tumor growth model calibration},
  author={Scheufele, Klaudius and Subramanian, Shashank and Mang, Andreas and Biros, George and Mehl, Miriam},
  journal={SIAM journal on scientific computing: a publication of the Society for Industrial and Applied Mathematics},
  volume={42},
  number={3},
  pages={B549},
  year={2020},
  publisher={NIH Public Access}
}
@inproceedings{gholami2019novel,
  title={A novel domain adaptation framework for medical image segmentation},
  author={Gholami, Amir and Subramanian, Shashank and Shenoy, Varun and Himthani, Naveen and Yue, Xiangyu and Zhao, Sicheng and Jin, Peter and Biros, George and Keutzer, Kurt},
  booktitle={Brainlesion: Glioma, Multiple Sclerosis, Stroke and Traumatic Brain Injuries: 4th International Workshop, BrainLes 2018, Held in Conjunction with MICCAI 2018, Granada, Spain, September 16, 2018, Revised Selected Papers, Part II 4},
  pages={289--298},
  year={2019},
  organization={Springer}
}
@article{avants2009advanced,
  title={Advanced normalization tools (ANTS)},
  author={Avants, Brian B and Tustison, Nick and Song, Gang and others},
  journal={Insight j},
  volume={2},
  number={365},
  pages={1--35},
  year={2009}
}
@article{swanson2000quantitative,
  title={A quantitative model for differential motility of gliomas in grey and white matter},
  author={Swanson, Kristin R and Alvord Jr, Ellsworth C and Murray, JD},
  journal={Cell proliferation},
  volume={33},
  number={5},
  pages={317--329},
  year={2000},
  publisher={Wiley Online Library}
}
@article{gooya2012glistr,
  title={GLISTR: glioma image segmentation and registration},
  author={Gooya, Ali and Pohl, Kilian M and Bilello, Michel and Cirillo, Luigi and Biros, George and Melhem, Elias R and Davatzikos, Christos},
  journal={IEEE transactions on medical imaging},
  volume={31},
  number={10},
  pages={1941--1954},
  year={2012},
  publisher={IEEE}
}
@article{zheng2014recovering,
  title={Recovering the source and initial value simultaneously in a parabolic equation},
  author={Zheng, Guang-Hui and Wei, Ting},
  journal={Inverse Problems},
  volume={30},
  number={6},
  pages={065013},
  year={2014},
  publisher={IOP Publishing}
}
@article{cheng2008quasi,
  title={A quasi Tikhonov regularization for a two-dimensional backward heat problem by a fundamental solution},
  author={Cheng, J and Liu, JJ},
  journal={Inverse Problems},
  volume={24},
  number={6},
  pages={065012},
  year={2008},
  publisher={IOP Publishing}
}
@book{ozisik2018inverse,
  title={Inverse heat transfer: fundamentals and applications},
  author={Ozisik, M Necat},
  year={2018},
  publisher={Routledge}
}
@inproceedings{ghafouri20233d,
  title={A 3D Inverse Solver for a Multi-species PDE Model of Glioblastoma Growth},
  author={Ghafouri, Ali and Biros, George},
  booktitle={International Workshop on Computational Mathematics Modeling in Cancer Analysis},
  pages={51--60},
  year={2023},
  organization={Springer}
}
@misc{hansen2022cma,
  title={CMA-ES/pycma: r3. 2.2},
  author={Hansen, N and Akimoto, Y and Baudis, P},
  year={2022},
  publisher={Zenodo}
}

\appendix
\section{Appendix}
\label{sec:appendix}
\begin{sidewaystable}
%  \begin{table}
  %\begin{adjustbox}{angle=90}
  \caption{Inversion results of parameters using synthetically generated data in 1D. In this table, we report the individual relative error for growth parameters for cases with noise levels. We report the relative error with regularization (Non-Reg) and without regularization (Reg). We also report the true values and the inverted values in the parenthesis in front of each relative error.}
  %\begin{table}
      \begin{adjustbox}{width=0.95\columnwidth,center}
        %\small \setlength{\tabcolsep}{1.85pt}      
        \begin{tabular}{cc|ccccccccccc}
              \toprule 
              Test-case & Noise (\%) & 
              $e_{\kappa}  (\frac{\kappa^{\rec}}{\kappa^{*}})$  &
              $e_{\rho} (\frac{\rho^{\rec}}{\rho^{*}}) $ & 
              $e_{\beta_0} (\frac{\beta_0^{\rec}}{\beta_0^{*}}) $ & 
              $e_{\alpha_0} (\frac{\alpha_0^{\rec}}{\alpha_0^{*}}) $ & 
              $e_{\gamma_0} (\frac{\gamma_0^{\rec}}{\gamma_0^{*}}) $ & 
              $e_{\delta_c} (\frac{\delta_c^{\rec}}{\delta_c^{*}}) $ &
              $e_{\delta_s} (\frac{\delta_s^{\rec}}{\delta_s^{*}}) $ &
              $e_{\ohyp} (\frac{\ohyp^{\rec}}{\ohyp^{*}}) $ &
              $e_{\oinv} (\frac{\oinv^{\rec}}{\oinv^{*}})$ & 
              $e_{\ith} (\frac{i^{\text{th},\rec}}{i^{\text{th},*}})$                                  \\ 
              \midrule
              Non-Reg & $5$ &${4.12e-01}$ $(\frac{8.47e-02}{6.00e-02})$ & 
              ${1.53e-01}$ $(\frac{1.61e+01}{1.90e+01})$ & 
              ${1.94e-01}$ $(\frac{6.45e+00}{8.00e+00})$ & 
              ${2.62e+00}$ $(\frac{7.24e+00}{2.00e+00})$ & 
              ${2.92e-02}$ $(\frac{1.03e+01}{1.00e+01})$ & 
              ${2.99e-01}$ $(\frac{9.11e+00}{1.30e+01})$ & 
              ${5.42e-01}$ $(\frac{1.37e+00}{3.00e+00})$ & 
              ${3.07e-01}$ $(\frac{3.92e-01}{3.00e-01})$ & 
              ${2.16e-01}$ $(\frac{3.92e-01}{5.00e-01})$ & 
              ${4.55e-01}$ $(\frac{4.37e-02}{3.00e-02})$ 
              \\ 
              %\rowcolor{gray!30}
              Reg &  $5$ & ${4.28e-01}$ $(\frac{8.57e-02}{6.00e-02})$ & 
              ${6.02e-02}$ $(\frac{1.79e+01}{1.90e+01})$ & 
              ${2.02e-01}$ $(\frac{6.38e+00}{8.00e+00})$ & 
              ${8.36e-01}$ $(\frac{3.27e-01}{2.00e+00})$ & 
              ${2.00e+00}$ $(\frac{3.00e+01}{1.00e+01})$ & 
              ${9.47e-02}$ $(\frac{1.42e+01}{1.30e+01})$ & 
              ${2.17e+00}$ $(\frac{9.51e+00}{3.00e+00})$ & 
              ${2.53e-01}$ $(\frac{3.76e-01}{3.00e-01})$ & 
              ${8.70e-02}$ $(\frac{5.43e-01}{5.00e-01})$ & 
              ${5.15e-01}$ $(\frac{4.54e-02}{3.00e-02})$ 
              \\ 
              \midrule
              Non-Reg & $10$ & ${5.46e-02}$ $(\frac{5.67e-02}{6.00e-02})$ & 
              ${2.56e-01}$ $(\frac{2.39e+01}{1.90e+01})$ & 
              ${1.37e+00}$ $(\frac{1.89e+01}{8.00e+00})$ & 
              ${9.48e-01}$ $(\frac{1.04e-01}{2.00e+00})$ & 
              ${1.54e-01}$ $(\frac{1.15e+01}{1.00e+01})$ & 
              ${7.25e-01}$ $(\frac{3.57e+00}{1.30e+01})$ & 
              ${6.67e-01}$ $(\frac{1.00e+00}{3.00e+00})$ & 
              ${1.33e+00}$ $(\frac{7.00e-01}{3.00e-01})$ & 
              ${9.66e-01}$ $(\frac{9.83e-01}{5.00e-01})$ & 
              ${3.86e-01}$ $(\frac{1.84e-02}{3.00e-02})$
              \\ 
              %\rowcolor{gray!30}
              Reg & $10$ & ${1.96e+00}$ $(\frac{1.77e-01}{6.00e-02})$ & 
              ${1.13e-01}$ $(\frac{2.11e+01}{1.90e+01})$ & 
              ${2.72e+00}$ $(\frac{2.98e+01}{8.00e+00})$ & 
              ${8.85e-01}$ $(\frac{2.31e-01}{2.00e+00})$ & 
              ${1.81e+00}$ $(\frac{2.81e+01}{1.00e+01})$ & 
              ${7.99e-01}$ $(\frac{2.34e+01}{1.30e+01})$ & 
              ${6.66e-01}$ $(\frac{1.00e+00}{3.00e+00})$ & 
              ${8.89e-01}$ $(\frac{3.33e-02}{3.00e-01})$ & 
              ${7.91e-02}$ $(\frac{4.60e-01}{5.00e-01})$ & 
              ${3.63e+00}$ $(\frac{1.39e-01}{3.00e-02})$ 
              \\
              \midrule
              Non-Reg &  $20$ & ${2.02e+00}$ $(\frac{1.81e-01}{6.00e-02})$ & 
              ${3.39e-01}$ $(\frac{1.26e+01}{1.90e+01})$ & 
              ${4.94e-01}$ $(\frac{4.05e+00}{8.00e+00})$ & 
              ${6.50e-01}$ $(\frac{7.01e-01}{2.00e+00})$ & 
              ${8.38e-01}$ $(\frac{1.84e+01}{1.00e+01})$ & 
              ${4.18e-01}$ $(\frac{7.57e+00}{1.30e+01})$ & 
              ${6.67e-01}$ $(\frac{1.00e+00}{3.00e+00})$ & 
              ${7.13e-01}$ $(\frac{5.14e-01}{3.00e-01})$ & 
              ${3.62e-01}$ $(\frac{6.81e-01}{5.00e-01})$ & 
              ${1.48e+00}$ $(\frac{7.43e-02}{3.00e-02})$
              \\
              %\rowcolor{gray!30}
              Reg & $20$ & ${1.17e+00}$ $(\frac{1.30e-01}{6.00e-02})$ & 
              ${1.03e-01}$ $(\frac{1.70e+01}{1.90e+01})$ & 
              ${3.96e-02}$ $(\frac{7.68e+00}{8.00e+00})$ & 
              ${9.50e-01}$ $(\frac{1.00e-01}{2.00e+00})$ & 
              ${2.00e+00}$ $(\frac{3.00e+01}{1.00e+01})$ & 
              ${4.57e-01}$ $(\frac{7.06e+00}{1.30e+01})$ & 
              ${6.67e-01}$ $(\frac{1.00e+00}{3.00e+00})$ & 
              ${5.84e-01}$ $(\frac{4.75e-01}{3.00e-01})$ & 
              ${7.49e-01}$ $(\frac{8.74e-01}{5.00e-01})$ & 
              ${3.26e-02}$ $(\frac{3.10e-02}{3.00e-02})$              
              \\
              \midrule
              Non-Reg &  $30$ & ${2.58e-01}$ $(\frac{7.55e-02}{6.00e-02})$ & 
              ${7.49e-02}$ $(\frac{1.76e+01}{1.90e+01})$ & 
              ${6.94e-02}$ $(\frac{8.55e+00}{8.00e+00})$ & 
              ${5.99e-01}$ $(\frac{8.02e-01}{2.00e+00})$ & 
              ${6.32e-01}$ $(\frac{3.68e+00}{1.00e+01})$ & 
              ${1.54e-03}$ $(\frac{1.30e+01}{1.30e+01})$ & 
              ${1.44e+00}$ $(\frac{7.31e+00}{3.00e+00})$ & 
              ${9.98e-01}$ $(\frac{5.99e-01}{3.00e-01})$ & 
              ${3.40e-01}$ $(\frac{6.70e-01}{5.00e-01})$ & 
              ${3.29e-01}$ $(\frac{2.01e-02}{3.00e-02})$
              \\
              %\rowcolor{gray!30}
              Reg &  $30$ & ${4.84e-01}$ $(\frac{8.90e-02}{6.00e-02})$ & 
              ${2.42e-02}$ $(\frac{1.85e+01}{1.90e+01})$ & 
              ${4.01e-01}$ $(\frac{1.12e+01}{8.00e+00})$ & 
              ${7.94e-01}$ $(\frac{4.13e-01}{2.00e+00})$ & 
              ${4.19e-01}$ $(\frac{5.81e+00}{1.00e+01})$ & 
              ${3.13e-01}$ $(\frac{8.93e+00}{1.30e+01})$ & 
              ${6.64e-01}$ $(\frac{1.01e+00}{3.00e+00})$ & 
              ${4.84e-01}$ $(\frac{4.45e-01}{3.00e-01})$ & 
              ${8.02e-01}$ $(\frac{9.01e-01}{5.00e-01})$ & 
              ${2.66e-01}$ $(\frac{2.20e-02}{3.00e-02})$                          
              \\                                                                                                  
              \toprule
          \end{tabular}      
      \end{adjustbox}
  \label{tab:1D_noise_err_params}
\end{sidewaystable}
%  \end{table} 
  
    %\begin{adjustbox}{angle=90}
    \begin{table}
        %\centering
        \caption{Inversion of synthetically generated data with $20\%$ noise level for $150$ parameters combinations with regularization. We report the average relative error for each individual parameter and also the average overall performance for all the parameters denoted as $e_{\mathbf{q}}$}
        \begin{adjustbox}{width=0.95\columnwidth,center}
          %\small \setlength{\tabcolsep}{1.85pt}      
          \begin{tabular*}{\textwidth}{c|cccccccccc|c}
                \toprule
                Regularization &                 
                $e_{\kappa}$  &
                $e_{\rho}$ & 
                $e_{\beta_0} $ & 
                $e_{\alpha_0} $ & 
                $e_{\gamma_0} $ & 
                $e_{\delta_c} $ &
                $e_{\delta_s}$ &
                $e_{\ohyp}$ &
                $e_{\oinv}$ & 
                $e_{\ith}$ &
                $e_{\mathbf{q}}$
                \\ 
                \midrule
                Non-Reg & $\num{1.42e+00}$ & $\num{3.11e-01}$ & $\num{5.80e+00}$ & $\num{1.68e+00}$ & $\num{1.53e+00}$ & $\num{2.55e+00}$ & $\num{2.23e+00}$ & $\num{9.24e-01}$ & $\num{4.14e-01}$ & $\num{1.12e+00}$ & $\num{1.80e+00}$ \\ 
                Reg & $\num{1.34e+00}$ & $\num{2.90e-01}$ & $\num{3.78e+00}$ & $\num{1.16e+00}$ & $\num{1.09e+00}$ & $\num{2.86e+00}$ & $\num{2.11e+00}$ & $\num{1.00e+00}$ & $\num{3.83e-01}$ & $\num{9.33e-01}$ & $\num{1.50e+00}$
                \\                                                                                             
                \bottomrule
            \end{tabular*}      
        \end{adjustbox}
    \label{tab:1D_noise_err_params_samples}
  \end{table}

  %\begin{adjustbox}{angle=90}
  \begin{table}
      \caption{Inversion results for individual model coefficients using synthetically generated data in 3D with known IC and brain anatomy. We report the relative error for each model coefficient and we report the inverted parameters and the true model coefficient. }
      \begin{adjustbox}{width=0.95\columnwidth,center}
        %\small \setlength{\tabcolsep}{1.85pt}      
          \begin{tabular}{cc|ccccccccccc}
              \toprule 
              Test-case & Noise (\%) 
              & $e_{\kappa}  (\frac{\kappa^{\rec}}{\kappa^{*}})$  &
              $e_{\rho} (\frac{\rho^{\rec}}{\rho^{*}}) $ & 
              $e_{\beta_0} (\frac{\beta_0^{\rec}}{\beta_0^{*}}) $ & 
              $e_{\alpha_0} (\frac{\alpha_0^{\rec}}{\alpha_0^{*}}) $ & 
              $e_{\gamma_0} (\frac{\gamma_0^{\rec}}{\gamma_0^{*}}) $ & 
              $e_{\delta_c} (\frac{\delta_c^{\rec}}{\delta_c^{*}}) $ &
              $e_{\delta_s} (\frac{\delta_s^{\rec}}{\delta_s^{*}}) $ &
              $e_{\ohyp} (\frac{\ohyp^{\rec}}{\ohyp^{*}}) $ &
              $e_{\oinv} (\frac{\oinv^{\rec}}{\oinv^{*}}) $ & 
              $e_{\ith} (\frac{i^{\text{th},\rec}}{i^{\text{th},*}})$ \\ 
              \midrule
              Non-Reg & 5 & $\num{2.58e-01}$ $(\frac{3.71e-03}{5.00e-03})$ & 
              $\num{2.28e-01}$ $(\frac{1.84e+01}{1.50e+01})$ & 
              $\num{2.51e+00}$ $(\frac{7.01e+00}{2.00e+00})$ & 
              $\num{1.99e+00}$ $(\frac{5.99e+00}{2.00e+00})$ & 
              $\num{4.89e-01}$ $(\frac{9.21e+00}{1.80e+01})$ & 
              $\num{6.13e+00}$ $(\frac{7.13e+00}{1.00e+00})$ & 
              $\num{6.44e-03}$ $(\frac{1.79e+01}{1.80e+01})$ & 
              $\num{1.11e+00}$ $(\frac{4.22e-01}{2.00e-01})$ & 
              $\num{4.93e-02}$ $(\frac{4.75e-01}{5.00e-01})$ & 
              $\num{1.83e-01}$ $(\frac{2.37e-02}{2.00e-02})$ \\ 
              %\rowcolor{gray!30}
              Reg & 5 & $\num{3.30e-01}$ $(\frac{6.65e-03}{5.00e-03})$ & 
              $\num{2.47e-03}$ $(\frac{1.50e+01}{1.50e+01})$ & 
              $\num{3.71e-03}$ $(\frac{2.01e+00}{2.00e+00})$ & 
              $\num{7.39e-02}$ $(\frac{2.15e+00}{2.00e+00})$ & 
              $\num{6.32e-02}$ $(\frac{1.91e+01}{1.80e+01})$ & 
              $\num{7.53e-01}$ $(\frac{1.75e+00}{1.00e+00})$ & 
              $\num{3.20e-01}$ $(\frac{1.22e+01}{1.80e+01})$ & 
              $\num{7.26e-01}$ $(\frac{3.45e-01}{2.00e-01})$ & 
              $\num{2.65e-01}$ $(\frac{6.32e-01}{5.00e-01})$ & 
              $\num{2.87e-02}$ $(\frac{2.06e-02}{2.00e-02})$ \\ 
              \midrule
              Non-Reg & 10 & $\num{3.60e-01}$ $(\frac{6.80e-03}{5.00e-03})$ & 
              $\num{5.47e-02}$ $(\frac{1.58e+01}{1.50e+01})$ & 
              $\num{3.66e-02}$ $(\frac{2.07e+00}{2.00e+00})$ & 
              $\num{1.56e-01}$ $(\frac{2.31e+00}{2.00e+00})$ & 
              $\num{2.21e-01}$ $(\frac{1.40e+01}{1.80e+01})$ & 
              $\num{6.14e+00}$ $(\frac{7.14e+00}{1.00e+00})$ & 
              $\num{1.10e-01}$ $(\frac{2.00e+01}{1.80e+01})$ & 
              $\num{1.82e-01}$ $(\frac{2.36e-01}{2.00e-01})$ & 
              $\num{2.97e-01}$ $(\frac{3.51e-01}{5.00e-01})$ & 
              $\num{1.53e-02}$ $(\frac{1.97e-02}{2.00e-02})$ \\ 
              %\rowcolor{gray!30}
              Reg & 10 & $\num{4.01e-01}$ $(\frac{7.01e-03}{5.00e-03})$ & 
              $\num{4.51e-02}$ $(\frac{1.57e+01}{1.50e+01})$ & 
              $\num{2.72e-02}$ $(\frac{2.05e+00}{2.00e+00})$ & 
              $\num{3.68e-01}$ $(\frac{2.74e+00}{2.00e+00})$ & 
              $\num{3.86e-02}$ $(\frac{1.73e+01}{1.80e+01})$ & 
              $\num{3.72e+00}$ $(\frac{4.72e+00}{1.00e+00})$ & 
              $\num{5.10e-01}$ $(\frac{8.82e+00}{1.80e+01})$ & 
              $\num{1.38e+00}$ $(\frac{4.76e-01}{2.00e-01})$ & 
              $\num{1.92e-01}$ $(\frac{5.96e-01}{5.00e-01})$ & 
              $\num{5.10e-03}$ $(\frac{1.99e-02}{2.00e-02})$ \\ 
              \midrule
              Non-Reg & 20 & $\num{4.76e-01}$ $(\frac{7.38e-03}{5.00e-03})$ & 
              $\num{2.58e-03}$ $(\frac{1.50e+01}{1.50e+01})$ & 
              $\num{2.81e-02}$ $(\frac{1.94e+00}{2.00e+00})$ & 
              $\num{8.59e-02}$ $(\frac{2.17e+00}{2.00e+00})$ & 
              $\num{3.59e-01}$ $(\frac{1.15e+01}{1.80e+01})$ & 
              $\num{6.95e+00}$ $(\frac{7.95e+00}{1.00e+00})$ & 
              $\num{2.10e-01}$ $(\frac{1.42e+01}{1.80e+01})$ & 
              $\num{1.04e+00}$ $(\frac{4.08e-01}{2.00e-01})$ & 
              $\num{2.10e-01}$ $(\frac{6.05e-01}{5.00e-01})$ & 
              $\num{6.10e-02}$ $(\frac{2.12e-02}{2.00e-02})$ \\ 
              %\rowcolor{gray!30}
              Reg & 20 & $\num{1.44e-01}$ $(\frac{5.72e-03}{5.00e-03})$ & 
              $\num{2.72e-02}$ $(\frac{1.54e+01}{1.50e+01})$ & 
              $\num{1.93e-02}$ $(\frac{1.96e+00}{2.00e+00})$ & 
              $\num{5.91e-02}$ $(\frac{2.12e+00}{2.00e+00})$ & 
              $\num{3.34e-01}$ $(\frac{1.20e+01}{1.80e+01})$ & 
              $\num{3.60e+00}$ $(\frac{4.60e+00}{1.00e+00})$ & 
              $\num{6.75e-02}$ $(\frac{1.68e+01}{1.80e+01})$ & 
              $\num{4.94e-01}$ $(\frac{2.99e-01}{2.00e-01})$ & 
              $\num{4.58e-02}$ $(\frac{4.77e-01}{5.00e-01})$ & 
              $\num{1.93e-02}$ $(\frac{2.04e-02}{2.00e-02})$ \\ 
              \midrule
              Non-Reg & 30 & $\num{1.09e+00}$ $(\frac{1.05e-02}{5.00e-03})$ & 
              $\num{3.53e-02}$ $(\frac{1.55e+01}{1.50e+01})$ & 
              $\num{3.58e-02}$ $(\frac{1.93e+00}{2.00e+00})$ & 
              $\num{1.13e+00}$ $(\frac{4.25e+00}{2.00e+00})$ & 
              $\num{7.54e-02}$ $(\frac{1.66e+01}{1.80e+01})$ & 
              $\num{6.23e-01}$ $(\frac{1.62e+00}{1.00e+00})$ & 
              $\num{7.56e-02}$ $(\frac{1.94e+01}{1.80e+01})$ & 
              $\num{4.67e-02}$ $(\frac{1.91e-01}{2.00e-01})$ & 
              $\num{5.45e-01}$ $(\frac{2.28e-01}{5.00e-01})$ & 
              $\num{1.00e-01}$ $(\frac{2.20e-02}{2.00e-02})$ \\ 
              %\rowcolor{gray!30}
              Reg & 30 & $\num{1.67e-03}$ $(\frac{4.99e-03}{5.00e-03})$ & 
              $\num{4.12e-02}$ $(\frac{1.56e+01}{1.50e+01})$ & 
              $\num{6.87e-02}$ $(\frac{1.86e+00}{2.00e+00})$ & 
              $\num{1.23e-01}$ $(\frac{1.75e+00}{2.00e+00})$ & 
              $\num{7.73e-02}$ $(\frac{1.94e+01}{1.80e+01})$ & 
              $\num{2.65e+00}$ $(\frac{3.65e+00}{1.00e+00})$ & 
              $\num{1.30e-02}$ $(\frac{1.82e+01}{1.80e+01})$ & 
              $\num{1.63e-01}$ $(\frac{2.33e-01}{2.00e-01})$ & 
              $\num{9.69e-02}$ $(\frac{4.52e-01}{5.00e-01})$ & 
              $\num{6.69e-02}$ $(\frac{2.13e-02}{2.00e-02})$ \\ 
              \toprule
          \end{tabular}      
      \end{adjustbox}
  %    \label{tab:notations}
  \label{tab:3D_noise_err_params}
  \end{table}
  %\end{adjustbox}  
    
    %\begin{adjustbox}{angle=90}
    \begin{table}
        \caption{Inversion results of model coefficients for synthetically generated data in 3D with unknown IC and brain anatomy. In this table, we report relative error for individual model coefficients and also we report the inverted model coefficients and the ground truth model coefficients in the parenthesis in from of errors. }
        \begin{adjustbox}{width=0.95\columnwidth,center}
         % \small \setlength{\tabcolsep}{1.85pt}      
            \begin{tabular}{c|ccccccccccccc}
                \toprule 
                Test-case 
                & $e_{\kappa}  (\frac{\kappa^{\rec}}{\kappa^{*}})$  &
                $e_{\rho} (\frac{\rho^{\rec}}{\rho^{*}}) $ & 
                $e_{\beta_0} (\frac{\beta_0^{\rec}}{\beta_0^{*}}) $ & 
                $e_{\alpha_0} (\frac{\alpha_0^{\rec}}{\alpha_0^{*}}) $ & 
                $e_{\gamma_0} (\frac{\gamma_0^{\rec}}{\gamma_0^{*}}) $ & 
                $e_{\delta_c} (\frac{\delta_c^{\rec}}{\delta_c^{*}}) $ &
                $e_{\delta_s} (\frac{\delta_s^{\rec}}{\delta_s^{*}}) $ &
                $e_{\ohyp} (\frac{\ohyp^{\rec}}{\ohyp^{*}}) $ &
                $e_{\oinv} (\frac{\oinv^{\rec}}{\oinv^{*}}) $ & 
                $e_{\ith} (\frac{i^{\text{th},\rec}}{i^{\text{th},*}})$          
                \\ 
                \midrule
  Case 1 & $\num{3.90e+01}$ $(\frac{2.00e-01}{5.00e-03})$ & 
  $\num{1.01e-01}$ $(\frac{1.35e+01}{1.50e+01})$ & 
  $\num{1.60e+00}$ $(\frac{5.20e+00}{2.00e+00})$ & 
  $\num{6.46e-01}$ $(\frac{7.09e-01}{2.00e+00})$ & 
  $\num{9.59e-02}$ $(\frac{1.63e+01}{1.80e+01})$ & 
  $\num{5.29e-02}$ $(\frac{1.05e+00}{1.00e+00})$ & 
  $\num{1.07e-01}$ $(\frac{1.99e+01}{1.80e+01})$ & 
  $\num{1.60e-01}$ $(\frac{2.32e-01}{2.00e-01})$ & 
  $\num{3.53e-01}$ $(\frac{3.24e-01}{5.00e-01})$ & 
  $\num{7.06e-01}$ $(\frac{5.89e-03}{2.00e-02})$ \\
  Case 2 &  $\num{5.41e+00}$ $(\frac{3.21e-02}{5.00e-03})$ & 
  $\num{1.29e-01}$ $(\frac{1.69e+01}{1.50e+01})$ & 
  $\num{4.23e+00}$ $(\frac{1.05e+01}{2.00e+00})$ & 
  $\num{6.37e-01}$ $(\frac{7.27e-01}{2.00e+00})$ & 
  $\num{3.34e-02}$ $(\frac{1.74e+01}{1.80e+01})$ & 
  $\num{1.15e-01}$ $(\frac{1.12e+00}{1.00e+00})$ & 
  $\num{6.38e-01}$ $(\frac{6.52e+00}{1.80e+01})$ & 
  $\num{1.44e+00}$ $(\frac{4.88e-01}{2.00e-01})$ & 
  $\num{9.18e-03}$ $(\frac{4.95e-01}{5.00e-01})$ & 
  $\num{9.50e-01}$ $(\frac{1.00e-03}{2.00e-02})$ \\
  Case 3 &  $\num{8.65e-01}$ $(\frac{9.32e-03}{5.00e-03})$ & 
  $\num{1.74e-01}$ $(\frac{1.76e+01}{1.50e+01})$ & 
  $\num{4.00e+00}$ $(\frac{1.00e+01}{2.00e+00})$ & 
  $\num{4.49e-01}$ $(\frac{1.10e+00}{2.00e+00})$ & 
  $\num{5.23e-02}$ $(\frac{1.71e+01}{1.80e+01})$ & 
  $\num{3.59e-01}$ $(\frac{1.36e+00}{1.00e+00})$ & 
  $\num{5.07e-01}$ $(\frac{8.87e+00}{1.80e+01})$ & 
  $\num{7.25e-01}$ $(\frac{3.45e-01}{2.00e-01})$ & 
  $\num{1.95e-01}$ $(\frac{4.03e-01}{5.00e-01})$ & 
  $\num{9.50e-01}$ $(\frac{1.00e-03}{2.00e-02})$ \\ 
                \toprule
            \end{tabular}      
        \end{adjustbox}
    %    \label{tab:notations}
    \label{tab:3D_noise_err_params_ic}
    \end{table}
  \begin{table}
      \caption{Inversion results for model coefficients using clinical data. For each patient, we report the inverted model coefficients using our inversion scheme.}
      \begin{adjustbox}{width=0.95\columnwidth,center}
      %  \small \setlength{\tabcolsep}{1.85pt}      
          \begin{tabular}{c|cccccccccccc}
              \toprule 
              Patient-ID
              & $\kappa $  &
              $ \rho $ & 
              $ \beta_0 $ & 
              $ \alpha_0 $ & 
              $ \gamma_0  $ & 
              $ \delta_c$ &
              $ \delta_s $ &
              $ \ohyp $ &
              $ \oinv $ & 
              $ \ith $   
              \\ 
              \midrule
              BraTS20\_Training\_039 & $\num{1.21e-03}$ & 
              $\num{1.50e+01}$ & 
              $\num{1.49e+01}$ & 
              $\num{9.98e+00}$ & 
              $\num{1.95e+01}$ & 
              $\num{1.42e+00}$ & 
              $\num{8.24e+00}$ & 
              $\num{2.64e-01}$ & 
              $\num{2.65e-01}$ & 
              $\num{1.12e-03}$ \\ 
              BraTS20\_Training\_042 & $\num{1.10e-03}$ & 
              $\num{1.88e+01}$ & 
              $\num{1.50e+01}$ & 
              $\num{1.31e+00}$ & 
              $\num{1.47e+01}$ & 
              $\num{2.49e+00}$ & 
              $\num{1.41e+01}$ & 
              $\num{3.15e-01}$ & 
              $\num{3.69e-01}$ & 
              $\num{2.15e-02}$ \\              
              \toprule
          \end{tabular}      
      \end{adjustbox}
  %    \label{tab:notations}
  \label{tab:3D_clinical_params}
  \end{table}

\end{document}